\documentclass[preprint,12pt]{elsarticle}

\usepackage[cmex10]{amsmath}
\usepackage{amssymb}
\usepackage{amsfonts}
\usepackage[T1]{fontenc}
\usepackage{textcomp}
\usepackage{lmodern}
\usepackage{array}
\usepackage{mathrsfs}
\usepackage{comment}
\usepackage{graphicx}
\usepackage{rotating}
\usepackage{bm}
\usepackage{url}
\usepackage{bm,amssymb,amsthm}
\usepackage{algorithm}
\usepackage{algorithmic}
\usepackage{multirow}
\usepackage{standalone}
\usepackage{blindtext}
\usepackage{float}
\usepackage{hyperref}
\usepackage{bm}
\usepackage{makecell}
\usepackage{array}
\usepackage{subscript}
\usepackage{textcomp}
\usepackage{caption}
\usepackage{setspace}
\usepackage{ctable}
\usepackage{enumitem}
\usepackage{listings}
\usepackage{xcolor}
\usepackage{subcaption}
\usepackage{array} 
\usepackage{booktabs} 

\definecolor{codegray}{rgb}{0.4,0.4,0.4}
\definecolor{backcolour}{rgb}{0.97,0.97,0.97}

\lstdefinestyle{journalpython}{
    language=Python,
    backgroundcolor=\color{backcolour},
    basicstyle=\ttfamily\scriptsize,  
    numbers=left,
    numberstyle=\tiny\color{codegray},
    stepnumber=1,
    numbersep=6pt,
    showspaces=false,
    showstringspaces=false,
    showtabs=false,
    tabsize=4,
    breaklines=true,
    breakatwhitespace=false,
    frame=single,                        
    rulecolor=\color{black},
    captionpos=b,
    aboveskip=0.5em,
    belowskip=0.5em
}

\newcommand{\bth}{\bm{\theta}}

\mathchardef\mhyphen="2D

\newcommand{\gr}[1]{\bm{#1}}

\newcommand{\q}{\gr{q}}
\newcommand{\Q}{\gr{Q}}

\newcommand{\A}{\gr{A}}
\newcommand{\B}{\gr{B}}

\newcommand{\I}{\gr{I}}

\newcommand{\C}{\gr{C}}
\newcommand{\0}{\gr{0}}
\newcommand{\Ss}{\gr{S}}
\newcommand{\D}{\gr{D}}

\newcommand{\R}{\gr{R}}
\newcommand{\N}{\gr{N}}
\renewcommand{\H}{\gr{H}}

\renewcommand{\S}{\gr{S}}
\newcommand{\G}{\gr{G}}

\newcommand{\U}{\gr{U}}

\newcommand{\V}{\gr{V}}

\newcommand{\Pre}{\gr{P}}

\newcommand{\n}{\gr{n}}

\newcommand{\x}{\gr{x}}

\newcommand{\y}{\gr{y}}

\newcommand{\df}{\overset{df.}{=}}

\newcommand{\sgens}[1]{\mathbb{R}^{#1}}
\newcommand{\sgen}[2]{\mathbb{R}^{#1\times #2}}

\newcommand{\sm}{\mathbb{R}^m}

\newcommand{\tr}{\mathbb{R}^{t\times r}}

\newcommand{\me}[1]{\mathbb{E}[{#1}]}

\newcommand{\ga}{\gamma}
\newcommand{\Si}{\gr{\Sigma}}

\newcommand{\la}{\lambda}
\newcommand{\Ga}{\gr{\Gamma}}

\newcommand{\no}{\noindent}
\newcommand{\pd}{\no\hspace{2em}{\itshape Proof: }}

\newtheorem{proposition}{Proposition}

\newtheorem{definition}{Definition}
\newtheorem{remark}{Remark}
\newtheorem{fact}{Fact}

\journal{NeuroImage}

\begin{document}

\begin{frontmatter}



\title{Multi-Source Neural Activity Indices for EEG/MEG Localization: A Two-Stage Spatial Filtering Framework and Extension to MNE-Python}
  
\author[PHYSUW]{Julia Jurkowska\fnref{affil,contr}}
\author[PHILUMK]{Joanna Dreszer}
\author[PHILUMK]{Monika Lewandowska}
\author[PHILUMK]{Krzysztof Tołpa}
\author[PHYSUMK,IAS]{Tomasz Piotrowski\corref{corr}\fnref{contr}}
\ead{tpiotrowski@umk.pl}
\cortext[corr]{Corresponding author}
\fntext[affil]{This work was initiated when Julia Jurkowska was a postgraduate student at the Faculty of Physics, University of Warsaw.}
\fntext[contr]{The first and last authors contributed equally to this work.}

\affiliation[PHYSUW]{organization={Faculty of Physics,
  University of Warsaw},
            addressline={Pasteura 5},
            city={Warsaw},
            postcode={02-093},
            country={Poland}}

\affiliation[PHILUMK]{organization={
    Faculty of Philosophy and Social Sciences,\\
    Nicolaus Copernicus University in Toruń},
            addressline={Fosa Staromiejska 1a},
            city={Toruń},
            postcode={87-100},
            country={Poland}}

\affiliation[PHYSUMK]{organization={
    Faculty of Physics, Astronomy and Informatics\\
    Nicolaus Copernicus University in Toruń},
            addressline={Grudziądzka 5},
            city={Toruń},
            postcode={87-100},
            country={Poland}}

\affiliation[IAS]{organization={
    Institute of Advanced Studies\\
    Nicolaus Copernicus University in Toruń, Poland},
            addressline={Wileńska 4},
            city={Toruń},
            postcode={87-100},
            country={Poland}}

\begin{abstract}
Accurate electroencephalography (EEG) and magnetoencephalography\linebreak (MEG) source localization and reconstruction are essential for understanding brain function, yet remain challenging because the underlying EEG/MEG inverse problem is inherently ill-posed. Spatial filtering (beamforming) approaches, such as linearly constrained minimum variance (LCMV) spatial filters, are widely used and well supported by existing analysis software.
In this work, we extend this framework by deriving a novel family of unbiased multi-source neural activity indices that form the localization stage of a two-stage spatial-filtering-based localization-reconstruction framework for the EEG/MEG inverse problem. In contrast to existing formulations, the proposed indices do not require knowledge of the target source covariance matrix, making them directly applicable in practical experimental settings. Their compact algebraic forms enable straightforward and numerically efficient implementation. The framework is validated on simulated EEG data and its applicability is illustrated through an example involving experimental EEG data from an oddball paradigm. To facilitate adoption, we provide a full open-source implementation extending MNE-Python, accompanied by a practical tutorial.
\end{abstract}




\end{frontmatter}



\section{Introduction}
Electroencephalography (EEG) and magnetoencephalography (MEG)\linebreak provide noninvasive access to neural dynamics by measuring the electromagnetic fields generated by synchronous post-synaptic currents in populations of neurons, using arrays of sensors positioned outside the brain \cite{Mosher1999,Michel2009}. Solving the associated inverse problem, i.e., inferring the location and dynamics of the underlying neural generators from sensor-level recordings, is therefore central to both scientific and medical applications \cite{Baillet2001}. Between the two modalities, EEG is far more widely used owing to its lower cost and portability compared with MEG. Thus, we focus below on EEG and present representative applications of solving the EEG inverse problem in both scientific and clinical contexts. It is important to underscore, however, that the methodological framework developed in this work applies to the MEG source localization problem without modification.

In neuroscience, solving the EEG inverse problem enables researchers to infer the origins of electrical activity by combining the millisecond-level temporal resolution of EEG with spatial accuracy that cannot be\linebreak achieved using sensor-level measurements alone. This approach has proven essential for identifying neural networks involved in cognitive processes such as attention, working memory, and motor control \cite{Michel2012, Mahjoory2017, Michel2019}. In clinical neuroscience, solving the EEG inverse problem has become a critical tool for the presurgical evaluation of patients with epilepsy (for a review, see \cite{Mierlo2020}). It facilitates the noninvasive localization of epileptogenic zones, thereby guiding surgical interventions with greater precision \cite{Brodbeck2010,Pelle2025}. Moreover, EEG inverse modeling is increasingly applied to the investigation of psychiatric and neurological disorders, including schizophrenia \cite{Canuet2011,Adamczyk2025}, ADHD \cite{Bluschke2018}, autism \cite{Wang2025}, and Alzheimer’s disease \cite{Aghajani2013}, where it aids in detecting altered neural connectivity and dysfunction within specific cortical networks.

One of the most widely used methods for solving the EEG/MEG inverse problem is the linearly constrained minimum variance (LCMV) spatial filter, introduced in array signal processing by Frost \cite{Frost1972}, and adapted to the EEG/MEG setup in, e.g., \cite{VanVeen1997, Gross2001, Greenblatt2005, Sekihara2008, Westner2022}. Its implementation is now a standard feature of all major EEG/MEG source-analysis toolboxes \cite{FieldTrip2011, Brainstorm2011, MNE2014, Jaiswal2019}. However, when the activity of multiple sources is correlated, the single-source formulation becomes suboptimal, as recognized since its introduction \cite{VanVeen1997}. Consequently, numerous spatial filtering approaches have been developed to explicitly address the correlated-source scenario \cite{Dalal2006,Brookes2007,Quraan2010,Diwakar2011,Moiseev2011,Piotrowski2019,Nunes2020,Kuznetsova2021}. In particular, some of these developments gave rise to source localization methods, i.e., focusing on spatial identification of active sources. Namely, source localization neural activity indices (NAIs) have been developed, e.g., based on the multi-source minimum variance (multi-source LCMV) filter \cite{Moiseev2011}, or its reduced-rank extensions \cite{Piotrowski2021}, based on the multi-source minimum-variance pseudo-unbiased reduced-rank (MV-PURE) framework \cite{Piotrowski2019, Yamada2006, Piotrowski2008, Piotrowski2009}. These neural activity indices generalize classical single-source NAIs \cite{VanVeen1997, Robinson1999} by operating in multi-dimensional signal subspaces, enabling robust localization in the presence of correlated sources.

However, despite their potential, multi-source localization methods for solving the EEG/MEG inverse problem have not yet gained wide acceptance within the neuroscientific community. This limited adoption can be attributed to the lack of intuitive use cases and accessible implementations of multi-source neural activity indices, which has restricted their practical applicability and left these approaches perceived more as intriguing alternatives than as integral components of commonly employed EEG/MEG source localization methodologies.

This work addresses this gap by providing a rigorous, step-by-step derivation of a family of novel, unbiased multi-source neural activity indices within an explicitly formulated two-stage localization-reconstruction framework for the EEG/MEG inverse problem. In contrast to prior approaches, in which these steps are often treated implicitly, we explicitly decouple the solution into two stages: (i) localization using the proposed family of multi-source neural activity indices, and (ii) reconstruction using spatial filters. Unlike the neural activity indices introduced in \cite{Piotrowski2021}, the proposed indices do not rely on the covariance matrix of the target sources, which is typically unavailable in practical experimental settings. Furthermore, they admit compact algebraic expressions, facilitating straightforward implementation. For the reconstruction stage, we use a standard LCMV beamformer as a reference spatial filter. This choice provides a well-established baseline; however, the framework is not restricted to LCMV and can accommodate alternative spatial filtering methods.

Furthermore, to support reproducibility and foster community adoption, we provide a full implementation of the proposed family of multi-source neural activity indices integrated into MNE-Python, one of the most widely used toolboxes for EEG/MEG analysis. The implementation is accompanied by an MNE-Python-style tutorial designed to facilitate uptake within the neuroscience community. The proposed indices are validated on simulated data, enabling controlled assessment of multi-source localization performance. We also demonstrate the practical utility of the framework through an application to experimental EEG data from a visual oddball paradigm.

The remainder of the paper is organized as follows. Section \ref{preli} introduces the notation, the EEG/MEG forward model, and the two-stage formulation of the solution to the EEG/MEG inverse problem. Section \ref{METHODS} presents the proposed family of multi-source neural activity indices and their integration into a spatial filtering-based localization-reconstruction framework. Section \ref{oddball} describes the MNE-Python implementation, including validation on simulated data, a tutorial, and an application to an EEG oddball dataset. Section~\ref{DISCUSSION} discusses methodological implications and future directions, with supplementary materials and proofs provided in the appendices. The ethics statement is provided in Section \ref{ethics}.

\section{Preliminaries} \label{preli}
\subsection{Forward Model}
We consider measurements of brain electrical activity by EEG/MEG sensors at a specified time interval. The random vector $\y\in\sm$ composed of measurements at a given time instant can be modeled as \cite{Mosher1999, VanVeen1997, Sekihara2008},
\begin{equation} \label{model}
\y=\H(\bth_0)\q_0+\n,
\end{equation}
where $\H(\bth_0)\in\sgen{m}{l_0}$ is the array response (lead-field) matrix of full column rank $l_0$, corresponding to sources at locations $\bth_0$, $\q_0\in\sgens{l_0}$ is a random vector representing dipole moments at $\bth_0$, and $\n\in\sm$ represents remaining brain activity along with noise recorded at the sensors.

We note that although model (\ref{model}) assumes constrained orientation of the sources, the results of this paper hold equally when orientations are unconstrained. In such a case, parametrization of a given source would include not only its position, but also the unit orientation vector. To simplify notation, we denote below $\H_0\df\H(\bth_0).$ 

We assume that $\q_0$ and $\n$ are uncorrelated zero-mean weakly stationary stochastic processes, so that
\begin{equation} \label{uncorr}
\me{\q_0\n^t}=\0\in\sgen{l_0}{m},
\end{equation}
and denote the positive definite covariance matrices of $\q_0$ and $\n$ as
\begin{equation} \label{QN}
\me{\q_0\q_0^t}\df\Q_0\succ\0,\quad\me{\n\n^t}\df\N\succ\0,
\end{equation}
respectively. Note that the assumption (\ref{uncorr}) implies that $\y$ is also zero-mean weakly stationary process with positive definite covariance matrix
\begin{equation} \label{R}
\me{\y\y^t}=\me{(\bm{H}_0\q_0+\n)(\bm{H}_0\q_0+\n)^t}=\H_0\Q_0\H_0^t+\N\df\R\succ\0.
\end{equation} 
We note that, while $\R$ can be estimated from the observed data and $\N$ is often estimable under suitable experimental conditions (e.g., from pre-stimulus activity in task-related experiments), the source covariance matrix $\Q_0$ is typically unknown.

\begin{remark} \label{ill}
It is important to emphasize a well-known fact that the EEG / MEG forward problem represented by model (\ref{model}) is inherently ill-posed, i.e., there are infinitely many alternative pairs $(\q,\n)$ which can result in the same measurements $\y.$ In other words, even under a  perfectly accurate forward model (\ref{model}), one cannot uniquely recover $\q_0$ without further constraints or priors, which is an unavoidable consequence of the underlying physics.
\end{remark}

\subsection{Two-Stage Solution to the EEG/MEG Inverse Problem} \label{TWO_STAGE}
We consider the following two‑stage solution to the EEG/MEG inverse problem:
\begin{enumerate}[label=\textcolor{red}{STAGE \arabic*},align=left,labelindent=1em]
\item Localize sources [determine $\bth_0$ in model (\ref{model})] using multi-source neural activity index, incorporating available neurophysiological temporal and spatial constraints; \label{localliza}
\item Reconstruct activity of sources at locations $\bth_0$ for a given time interval using spatial filter. \label{reconstra}
\end{enumerate}
Such a two-stage approach to the EEG/MEG inverse problem is implicit in most spatial filtering-based methods employing classical single-source neural activity indices \cite{VanVeen1997, Robinson1999}. The proposed solution makes this strategy explicit by formulating it as a two-stage localization-reconstruction framework and by introducing novel family of multi-source neural activity indices for the localization stage. This enables the joint evaluation of candidate source sets and substantially reduces the number of sources to be considered in the subsequent reconstruction stage. Further details are provided in Section \ref{METHODS}.

\subsection{Unbiased Neural Activity Indices}
The goal of multi-source neural activity indices is to identify the source location $\bth_0$ in model (\ref{model}) in the presence of the ill-posedness of the EEG/MEG forward problem (see Remark \ref{ill}). We therefore adopt the following definition, requiring the indices to be unbiased, i.e., $\bth_0$ belongs to the set of maximizers.

\begin{definition} \label{nai}
Let $\bm{\Theta}$ be the power set of feasible locations of sources and let $\bth\in\bm{\Theta}.$
We call a function $f:\bm{\Theta}\to\mathbb{R}_+$ a neural activity index (NAI). We say that a NAI $f$ is unbiased if~\cite{Sekihara2008,Moiseev2011,Piotrowski2021}
\begin{equation} \label{unbiasedness}
\bth_0\in\arg\max_{\bth\in\bm{\Theta}}f(\bth).
\end{equation}
\end{definition}

\section{Methods} \label{METHODS}
We first establish the necessary notation.

\begin{table}[h]
\centering
\caption{Matrices $\bm{G}(\bth)$ and $\bm{S}(\bth)$, together with their counterparts obtained via congruence transformation with $\bm{Q}(\bth)^{1/2}$.}
\label{tab:GST}
{\small
\renewcommand{\arraystretch}{1.8}
\setlength{\arrayrulewidth}{1pt}
\[
\begin{array}{|l|l|}
\hline
\bm{G(\bth)} \df \bm{H(\bth)}^t \bm{N}^{-1}\bm{H(\bth)}\in\sgen{l}{l}
&
\tilde{\G}(\bth)\df \Q(\bth)^{1/2}\G(\bth)\Q(\bth)^{1/2}\in\sgen{l}{l}
\\
\hline
\bm{S(\bth)} \df \bm{H(\bth)}^t \bm{R}^{-1}\bm{H(\bth)}\in\sgen{l}{l}
&
\tilde{\S}(\bth)\df \Q(\bth)^{1/2}\S(\bth)\Q(\bth)^{1/2}\in\sgen{l}{l}
\\
\hline
\end{array}
\]
}
\end{table}

\noindent In the above,
\begin{equation} \label{tQ}
\Q(\bth)\df \S(\bth)^{-1}-\G(\bth)^{-1}\in\sgen{l}{l},
\end{equation}
$\H(\bth)\in\sgen{m}{l}$ is the array response (lead-field) matrix of full column rank $l$, corresponding to sources at locations $\bth$, with $\R$ and $\N$ defined in (\ref{R}) and (\ref{QN}), respectively. We note that $\bm{G(\bth)}$ and $\bm{S(\bth)}$ are positive-definite for all $\bth\in\bm{\Theta}$ from Fact \ref{Horn772} in \ref{kru}.
Further, it is natural to assume that the activity of sources are not linearly dependent (almost surely), so that
\begin{equation} \label{NATURA}
\forall\bth\in\bm{\Theta}\quad \Q(\bth)\succ\0.
\end{equation}
From Fact \ref{Horn772} in \ref{kru}, this implies that also $\tilde{\G}(\bth)$ and $\tilde{\S}(\bth)$ are positive definite for all $\bth\in\bm{\Theta}.$ We note that $$\tilde{\G}(\bth)^{-1}=\Q(\bth)^{-1/2}\G(\bth)^{-1}\Q(\bth)^{-1/2}$$ and $$\tilde{\S}(\bth)^{-1}=\Q(\bth)^{-1/2}\S(\bth)^{-1}\Q(\bth)^{-1/2}.$$ For further use, we also define
\begin{equation} \label{G_Filt}
\bm{G}_0\df \bm{H}_0^t\bm{N}^{-1}\bm{H}_0\succ\0,
\end{equation}
and
\begin{equation} \label{S_Filt}
\bm{S}_0\df \bm{H}_0^t\bm{R}^{-1}\bm{H}_0\succ\0.
\end{equation}

\noindent Unless stated otherwise, we assume that the eigenvalues of a diagonalizable square matrix are arranged in decreasing order.

\subsection{Novel $MAI_{MVP}$ Neural Activity Indices} \label{proposed_NEW}
We introduce the following family of neural activity indices as solutions to the localization problem in \ref{localliza}, derived from the MV-PURE family of filters \cite{Piotrowski2019} (see discussion below Proposition \ref{FULLFORM2025}). Unlike the ``first-generation'' reduced-rank-capable neural activity indices from \cite{Piotrowski2021}, the proposed $MAI_{MVP}$ indices depend solely on the known component matrices $\G(\bth)$ and $\S(\bth)$, i.e., they do not require knowledge of the (typically unknown) $\Q_0$.
\pagebreak

\begin{definition} \label{DFMVP}
Consider $l$ sources of brain's electrical activity and let $1\leq r\leq l_0.$ Using notation introduced in Table \ref{tab:GST}, we define $MAI_{MVP}$ multi-source neural activity indices as
\begin{multline} \label{MAIMVP_small}
  MAI_{MVP}(\bth,r)\df tr\{\tilde{\G}(\bth)\tilde{\S}(\bth)^{-1}\}-l=\\tr\{\G(\bth)\S(\bth)^{-1}\}-l=
  \sum_{i=1}^{l}\lambda_i(\G(\bth)\S(\bth)^{-1})-l,\quad l\leq r,
\end{multline}
and
\begin{equation} \label{MAIMVP}
MAI_{MVP}(\bth,r)\df tr\{\tilde{\G}(\bth)\tilde{\S}(\bth)^{-1}\Pre^{(r)}_{\tilde{\S}(\bth)}\}-r,\quad l>r,
\end{equation}
where $\Pre^{(r)}_{\tilde{\S}(\bth)}$ is the orthogonal projection matrix onto subspace spanned by eigenvectors corresponding to the $r$ largest eigenvalues of $\tilde{\S}(\bth).$
\end{definition}

Proposition \ref{FULLFORM2025} establishes the unbiasedness of the $MAI_{MVP}$ neural activity indices and provides simple algebraic forms that we recommend for implementation. These forms provide an alternative to the MV-PURE-based interpretation of the proposed indices and are discussed below the proposition.

\begin{proposition} \label{FULLFORM2025}
Consider $l$ sources of brain's electrical activity and let $1\leq r\leq l_0.$ Then: 
\begin{enumerate}
\item $MAI_{MVP}$ neural activity indices are unbiased for any rank $r$, and the following equality holds: \label{MVP_unbiased}
\begin{equation} \label{th1}
MAI_{MVP}(\bth_0,r)=\sum_{i=1}^{r}\lambda_{i}(\R\N^{-1})-r,\quad r=1,\dots,l_0.
\end{equation}
\item For a given rank $r$ and $l>r$, the following equality holds: \label{MVP_MAX}
\begin{equation} \label{MAIMVPEXT}
MAI_{MVP}(\bth,r)=\sum_{i=1}^{r}\lambda_i(\G(\bth)\S(\bth)^{-1})-r.
\end{equation}
\end{enumerate}
\end{proposition}
\pd See \ref{PD_FULLFORM2025}.
\pagebreak

A few remarks on the obtained family of the reduced-rank multi-source neural activity indices $MAI_{MVP}$ are in place here. First, from (\ref{th1}) and (\ref{MAIMVPEXT}) we obtain for $\bth_0$ that
\begin{equation} \label{eigs_RN}
\sum_{i=1}^{r}\lambda_{i}(\R\N^{-1})
=
\sum_{i=1}^{r}\lambda_i(\G(\bth_0)\S(\bth_0)^{-1}),\quad r=1,\dots,l_0.
\end{equation}
Thus, the corresponding eigenvalues coincide; it can be also verified that the number of active sources $l_0$ is equal to the number of eigenvalues of $\R\N^{-1}$ that are strictly greater than one. This follows from the decomposition $\R=\H_0\Q_0\H_0^t+\N$ [see (\ref{R})], which implies that $\R\N^{-1}$ can be interpreted as a sensor-level signal-to-noise ratio. Under this interpretation, the number of eigenvalues strictly greater than one determines the dimension of the signal subspace.\footnote{Accordingly, we will also refer below to $l_0$ as the number of directions in the source space.} From this perspective,
$$
\G(\bth_0)\S(\bth_0)^{-1}
= \H_0^t\N^{-1}\H_0(\H_0^t\R^{-1}\H_0)^{-1}
$$
can be interpreted as a source-level signal-to-noise ratio, linked to the source space via the lead-field $\H_0$. As Proposition \ref{FULLFORM2025} shows, maximizing this ratio with respect to the source locations $\bth$ yields unbiased multi-source neural activity indices in (\ref{MAIMVPEXT}) for any principal $r$-dimensional subspace of the signal space, where $r = 1, \dots, l_0$.

We also note that the form of the $MAI_{MVP}$ neural activity indices in \ref{MAIMVP} admits an interpretation as source-level signal-to-noise-ratios at the outputs of the multi-source MV-PURE filters, which were shown in \cite{Piotrowski2019} to be highly robust to source correlation. Since multi-source spatial filtering has generally been shown to mitigate the effects of source correlation \cite{Dalal2006,Brookes2007,Quraan2010,Diwakar2011,Moiseev2011,Piotrowski2019,Nunes2020,Kuznetsova2021}, maximizing these output signal-to-noise ratios suggests that the resulting multi-source neural activity indices are less susceptible to signal cancellation caused by correlated sources.

Finally, we note that the $MAI_{MVP}$ indices provide a new interpretation of multi-source neural activity indices previously introduced in the literature:
\begin{itemize}
\item When $r=l$, so that no rank constraint is imposed, $MAI_{MVP}$ coincides with the $MAI$ neural activity index introduced in \cite{Moiseev2011}. Thus, the present work can also be viewed as extending $MAI$ to a broader family of unbiased multi-source neural activity indices.
\item The formulation of $MAI_{MVP}$ given in Proposition \ref{FULLFORM2025} sheds also new light on the findings of \cite{Moiseev2013}, where the eigenvalues of a related source-level signal-to-noise ratio were shown to have stationary points at $\bth_0.$ This property was referred to as unbiasedness in a broad sense in \cite{Moiseev2013}.
\end{itemize}

\subsection{Two-Stage Spatial-Filtering-Based Localization-Reconstruction Pipeline} \label{MNE_IMP}
Based on Section \ref{proposed_NEW}, and in accordance with the architecture introduced in Section \ref{TWO_STAGE}, we propose the following two-stage spatial-filtering-based localization-reconstruction framework for the EEG/MEG inverse problem:

\begin{enumerate}[label=\textcolor{blue}{STAGE \arabic*},align=left,labelindent=1em]
\item Proposition \ref{FULLFORM2025}: Localize the active sources $\bth_0$ using $MAI_{MVP}$ neural activity indices. To overcome the numerically intractable number of combinations needed for brute-force evaluation of multi-source indices, we note that $MAI_{MVP}$ are amenable to the efficient iterative implementation scheme introduced in \cite{Moiseev2011} and also used heuristically in \cite{Piotrowski2021}. For $MAI_{MVP}$, this iterative source search scheme is particularly well suited from a computational perspective (see \ref{ISSA} for details). It also possesses a desirable property of identifying sources sequentially in order of decreasing strength, enabling full utilization of the unbiasedness property of the proposed neural activity indices across all ranks. Specifically, at the first iteration all indices are rank-1 and identical, taking the form  (\ref{MAIMVP_small}), which leads to the detection of the strongest source. At the second iteration, the rank-1 index transitions to the form (\ref{MAIMVP}), while the remaining higher-rank indices retain the form (\ref{MAIMVP_small}). At the third iteration, the rank-1 and rank-2 indices take the form (\ref{MAIMVP}), whereas the higher-rank indices still follow the form (\ref{MAIMVP_small}). This pattern continues progressively, with an increasing number of indices adopting the form (\ref{MAIMVP}), until in the final iteration all indices-except the full-rank one-are of the form (\ref{MAIMVP}), thereby operating independently to identify distinct sources.\footnote{In practice, the form (\ref{MAIMVPEXT}) is used in place of (\ref{MAIMVP}).} Thus, this procedure introduces, by design, an increasing level of redundancy in the source search (it produces at most $l_0(l_0+1)/2$ distinct candidate sources in $l_0$ iterations). This redundancy promotes robust localization of weaker sources, while ensuring that the strongest sources are reliably identified in the initial iterations.

We also note that, in practice, preprocessing in EEG/MEG source localization replaces the (post-stimulus) signal covariance matrix $\R$ with its reduced-rank version $\tilde{\R}$ such that the rank of $\tilde{\R}$ is an effective dimension $d_{eff}$, corresponding to the subspace representable by the forward model. All subsequent quantities are therefore implicitly defined within this subspace. More precisely, for a set of candidate sources $\boldsymbol{\theta}$ with lead-field $\mathbf{H}(\boldsymbol{\theta}) \in \mathbb{R}^{m \times \ell}$ of full column rank $\ell \le d_{eff}$, active sources give rise to lead-field components that lie within the signal subspace, i.e., $\mathcal{R}(\mathbf{H}(\boldsymbol{\theta})) \subseteq \mathcal{R}(\tilde{\R})$. As a result, the matrix $\H(\boldsymbol{\theta})^t \tilde{\R} \H(\boldsymbol{\theta})$
is positive definite. Indeed, for any nonzero vector $\mathbf{y}$,
\[
\mathbf{y}^t \H(\boldsymbol{\theta})^t \tilde{\R} \H(\boldsymbol{\theta}) \mathbf{y}
= (\H(\boldsymbol{\theta}) \mathbf{y})^t \tilde{\R} (\H(\boldsymbol{\theta}) \mathbf{y}) > 0,
\]
since $\H(\boldsymbol{\theta}) \mathbf{y} \in \mathcal{R}(\tilde{\R})$ and $\tilde{\R}$ is positive definite on its range. This property ensures that components of the proposed neural activity indices involving the preprocessed $\tilde{\R}$ (such as $\tilde{\S}(\bth)\df\H(\bth)^t\tilde{\R}^{-1}\H(\bth)$) are well-defined and positive definite for active sources. For components that additionally depend on the noise covariance $\N$, it is natural to assume that $\tilde{\N}$ is non-singular on the same subspace, as it undergoes the same preprocessing and rank reduction as $\tilde{\R}$, e.g., when estimated from pre-stimulus data.

In contrast, inactive sources are associated with directions that are weakly represented in the data, which may lead to rank-deficient or poorly conditioned matrices such as $\tilde{\S}(\bth).$ This observation is consistent with interpreting the proposed neural activity indices as identifying directions that maximize the source-level signal-to-noise ratio, as discussed below Proposition \ref{FULLFORM2025}. \label{localize!}

\item  With the source locations $\bth_0$ established in the previous stage, we reconstruct the activity of sources at $\bth_0$ using the LCMV spatial-filtering framework implemented in MNE-Python. Although the proposed framework could naturally be extended to employ multi-source spatial filters in this stage, we adopt LCMV due to its widespread use and highly optimized implementation in MNE-Python. Importantly, the multi-source localization stage reduces the candidate source space from several thousand locations to a small set of sources, substantially improving the clarity and interpretability of the reconstruction stage. Because the neural activity indices derived in Section~\ref{proposed_NEW} are integrated with MNE-Python, the proposed two-stage localization-reconstruction framework directly extends its existing functionality. \label{reconstruct!}
\end{enumerate}

\section{Extending MNE-Python: Simulation Study and Application to Experimental Data with Tutorial} \label{oddball}

\subsection{Simulation Study} \label{sims}
To quantitatively evaluate the spatial accuracy of the proposed source localization approach in \ref{localize!}, we conducted validation experiments using simulated EEG data with known ground truth. Simulations enable controlled assessment of localization performance by providing direct access to the true cortical sources and their time courses. The simulation framework was implemented in Python using the MNE-Python library; the full code is available in the project repository
    \vspace{\baselineskip}
    \\\url{https://github.com/julia-jurkowska/mvpure-tools},
    \vspace{\baselineskip}

    \noindent along with a tutorial showcasing the use of the code
    \vspace{\baselineskip}
    \\\url{https://julia-jurkowska.github.io/mvpure-tools}.
    \vspace{\baselineskip}

\noindent The provided code and tutorial allow users to freely explore the full set of simulation parameters summarized in Table~\ref{tab:constant-parameters-simulation} in Section~\ref{sim_sup}.

\subsubsection{Simulation framework}
In the simulation framework, we aimed to replicate the experimental conditions described in Section \ref{exp_data}, based on EEG recordings from a visual oddball paradigm. To this end, source-level neural activity was generated and projected to sensor space using the same forward model as for the experimental data. The simulation comprised three stages: (1) selection of cortical source locations, (2) generation of source time series, and (3) projection to EEG sensor space with realistic noise.

\paragraph{Source Space and Anatomical Regions}
Cortical sources were defined within a surface-based source space derived from the forward model that is used for the experimental data in Section \ref{exp_data}. Anatomical regions were defined using the cortical parcellation provided by FreeSurfer with annotations from the Desikan-Killiany Atlas \cite{Desikan2006}. Sources were drawn from two sets of regions: a set of background labels covering a broad range of cortical areas, and its subset with stimulus-active labels, see Table \ref{tab:used_labels_desc} in \ref{sim_sup}.

\paragraph{Generation of Source Activity}
Source time series were simulated using multivariate autoregressive (MVAR) models with two types of activity:
\vspace{-.25\baselineskip}
\begin{itemize}
\item \textit{Background activity}, representing ongoing neural dynamics unrelated to
stimulus processing, generated across all selected sources with moderate cross-source
coupling, then temporally smoothed using a Gaussian filter and scaled to a physiologically
plausible amplitude range.
\vspace{-.25\baselineskip}
\item \textit{Stimulus-evoked activity}, generated within a predefined post-stimulus window ($50$-$200$~ms after stimulus onset) using a low-dimensional latent process. Latent signals were simulated with an MVAR model and projected onto the stimulus-active sources using a column-normalised randomly generated model matrix, ensuring consistent spatial activation pattern across epochs while still producing trial-to-trial amplitude variability.
\end{itemize}
\vspace{-.2\baselineskip}
In addition to the stochastic MVAR activity, a synthetic event-related potential (ERP) component (a P1-N1-P2 complex modeled as a sum of Gaussian functions) was added only to the stimulus-related sources. To increase realism, the ERP amplitude and latency were varied across epochs through random jitter. Furthermore, the ERP amplitude was modulated according to the anatomical region of each source, with occipital and inferior parietal regions receiving higher weights, reflecting the known topography of visual evoked responses.

In summary, the pre-stimulus background activity was used to estimate the noise covariance matrix $\N$, whereas the post-stimulus background activity with added stimulus-evoked activity and synthetic ERP component was used to estimate the signal covariance matrix $\R.$ We note that, in the latter case, we also evaluate the robustness of the proposed indices to mild violations of the zero-mean assumption due to the ERP component.

After generating the post-stimulus activity, its amplitude was scaled to achieve a predefined signal-to-noise ratio (SNR) at the source level, defined as
\begin{equation} \label{SNR}
\mathrm{SNR}_{\mathrm{dB}} = 10\log_{10}\left(\frac{P_{\mathrm{signal}}}{P_{\mathrm{bkg}}}\right),
\end{equation}
where $P_{\mathrm{signal}}$ is the power of the post-stimulus signal and $P_{\mathrm{bkg}}$ is the power of the background signal in the post-stimulus window, both averaged across $N$ simulated trials.\footnote{Alternative definitions of SNR may also be considered, for example by varying the number of trials $N$ used for estimation. The provided source code exposes all relevant parameters to the user for experimentation.}

\paragraph{Projection to Sensor Space}
The simulated cortical sources were projected to the EEG sensors using the subject-specific lead-field matrix employed for the experimental data in Section \ref{exp_data}. Two types of sensor-level noise were then added: low-frequency drift (generated by temporally smoothing Gaussian noise) and additive white Gaussian noise scaled relative to the signal amplitude. The resulting signals were organized into epochs and processed with band-pass filtering ($1$-$45$~Hz, IIR) and average-referenced.

\subsubsection{Validation Procedure}
Noise and signal covariance matrices were estimated from the pre-stimulus interval ($[-200,\ 0]$~ms) and the post-stimulus window ($[50,\ 200]$~ms), respectively. The localization stage described in \ref{localize!} using $MAI_{MVP}$ was then applied to the simulated EEG data, and the resulting source localizations were compared with the known ground truth. For each simulation run, the iterative localization algorithm followed the recommendation given in \ref{localize!} and produced at most $l_0(l_0+1)/2$ distinct candidate sources, which may overlap if $MAI_{MVP}$ of different ranks points to the same source. Here, $l_0$ denotes the total number of simulated active sources, which is also the number of directions in the source space (see discussion below Proposition~\ref{FULLFORM2025}).

To determine whether a source was correctly localized, spatial localization accuracy was quantified using a distance-based error metric representing the average spatial deviation between true and estimated sources. To this end, the true and estimated vertices were converted to their corresponding three-dimensional cortical coordinates using the source space geometry. A pairwise Euclidean distance matrix between true and estimated source locations was then computed, and optimal one-to-one matching between the two sets was performed using the Hungarian assignment algorithm to minimize the total localization error.

As a baseline comparison, the LCMV beamformer was applied to the same simulated data using the Neural Activity Index (NAI) as the standard localizer. The $l_0$ sources with the highest NAI values across the post-stimulus window were selected as the LCMV localization estimate.

The experiment was designed to assess how robust source localization is under different levels of signal-to-noise ratio (SNR) and varying source configurations. Two scenarios were considered: (1) two active regions, each consisting of two vertices separated by at most 0.75 cm (Euclidean distance); and (2) six active regions, with one vertex per region (see Table \ref{tab:used_labels_desc} in \ref{sim_sup}). Each scenario was tested across three SNR levels: 1 dB, 3 dB, 5 dB. For every combination of scenario and SNR, $N=100$ independent simulation runs were performed, with source locations randomly resampled for each run and region. All simulation parameters are summarized in Table \ref{tab:constant-parameters-simulation} in \ref{sim_sup}.

\begin{figure}[h!]
\centering
\includegraphics[width=1\linewidth]{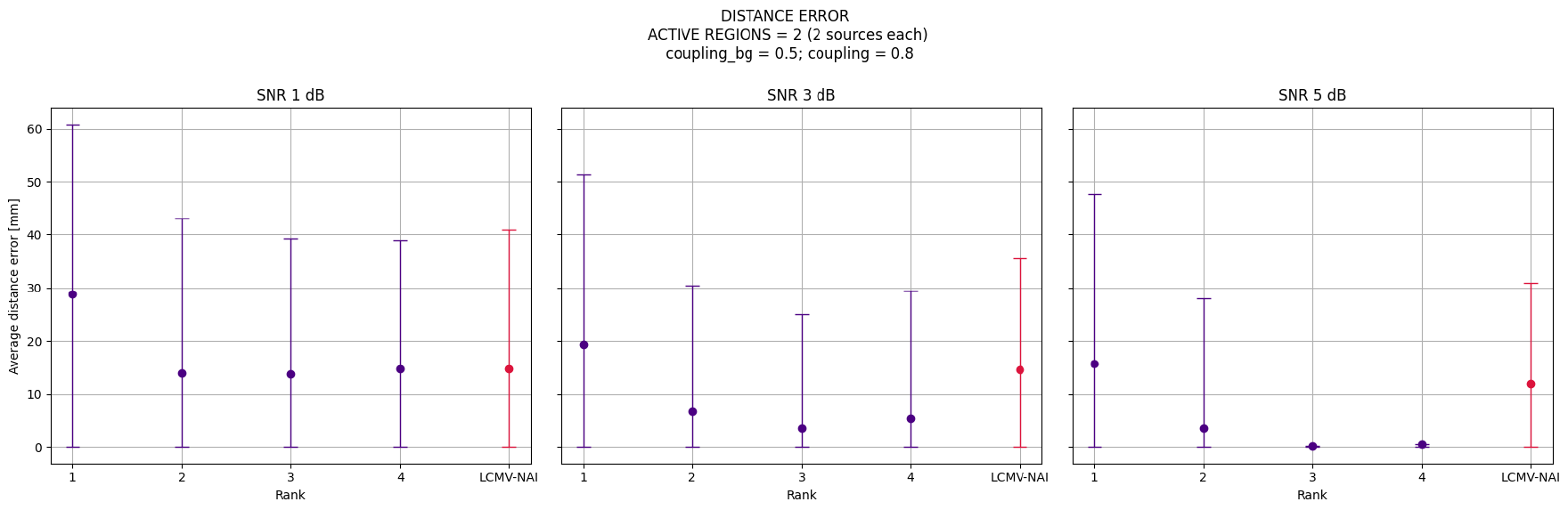}
\caption{Scenario (1): Average localization error [mm] across ranks of $MAI_{MVP}$ and $LCMV-NAI$, evaluated over a range of SNR values (Definition~\ref{SNR}), for two pairs of closely positioned sources (regions listed in Simulation~1, Table~\ref{tab:used_labels_desc} in \ref{sim_sup}).} \label{2x2}
\end{figure}

\begin{figure}[h!]
\centering
\includegraphics[width=1\linewidth]{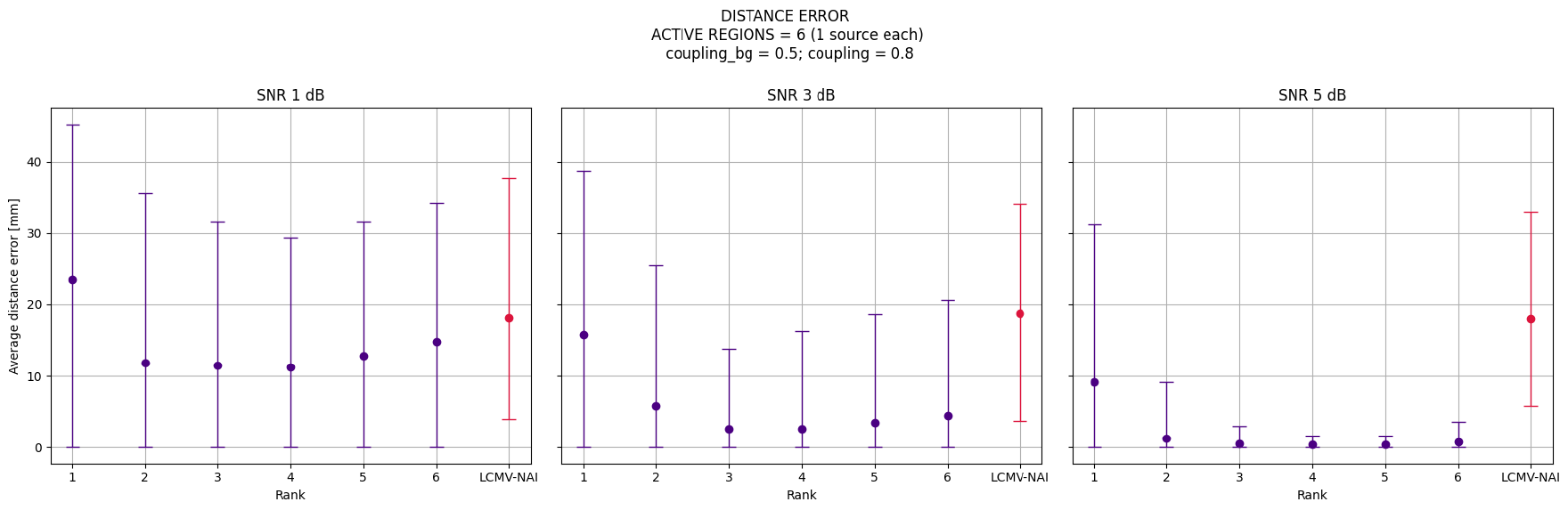}
\caption{Scenario (2): Average localization error [mm] across ranks of $MAI_{MVP}$ and $LCMV-NAI$, evaluated over a range of SNR values (Definition~\ref{SNR}), for six active sources (one per region; regions listed in Simulation~2, Table~\ref{tab:used_labels_desc} in \ref{sim_sup}).} \label{6x1}
\end{figure}

The performance of the $MAI_{MVP}$ indices is largely stable across ranks, with mid-range ranks yielding the best spatial accuracy in the considered setup, while the rank-one index consistently shows the poorest performance. However, for different data and experimental settings, different ranks may yield optimal localization performance. Therefore, we recommend following the procedure described in \ref{localize!}, i.e., evaluating $MAI_{MVP}$ across all ranks during the localization stage to fully exploit the unbiasedness property of the indices, which, by Proposition \ref{FULLFORM2025}, holds for all admissible ranks.

For the LCMV beamformer, the single-source $NAI$ is not an unbiased index, as reflected in the results above. In practice, LCMV-based localization is often supplemented by scanning the source space and identifying local maxima of the $NAI$. However, in the presence of strongly correlated sources, this strategy is of limited utility, as source cancellation effects lead to attenuated $NAI$ values. Here, we avoid such peak-selection heuristics to highlight the advantages of the proposed localization procedure in \ref{localize!}, which exploits redundant localization using a family of unbiased multi-source neural activity indices. 

\subsection{Application to Experimental Data with Tutorial} \label{exp_data}
As a key practical contribution of this work, we also provide a complete implementation of the methodology introduced in Section \ref{METHODS} to experimental data, available at:
\vspace{\baselineskip}
\\\url{https://github.com/julia-jurkowska/mvpure-tools}.
\vspace{\baselineskip}

\noindent The implementation is provided as an extension that integrates with MNE-Python and builds on its established API and tools (including covariance matrix regularization, lead-field preprocessing, forward/inverse model handling, and related utilities). To accelerate adoption, we include a step-by-step tutorial that walks through a representative application of the processing pipeline introduced in Section \ref{MNE_IMP}. The fully implemented tutorial (including preprocessing and link to dataset) is available at
\vspace{\baselineskip}
\\\url{https://julia-jurkowska.github.io/mvpure-tools}.
\vspace{\baselineskip}

The tutorial is provided on EEG data with surface lead-field model of fixed orientation normal to the cortex surface, with 5124 candidate sources.

\noindent Below, we provide a brief description of the experimental protocol and dataset, as well as the key excerpts of the tutorial which implement the source localization and reconstruction framework introduced in Section \ref{MNE_IMP}.

\subsubsection{Experimental Protocol and Dataset} \label{EP}
The dataset comprises data from high-density 128-channel EEG recordings acquired during a visual 3-stimulus oddball paradigm from a single right-handed, 27-year-old male participant, drawn from the publicly available Nencki-Symfonia dataset \cite{Dzianok2022} (see Ethics statement in Section \ref{ethics}). The participant completed a visual three-stimulus oddball paradigm, where three types of 200-ms visual stimuli (rare target - \th{}, deviant (distractor) - \TH, and frequent standard - p) were presented in the center of a screen in a pseudorandom order, such that two identical rare stimuli did not occur consecutively. Stimuli were presented using Presentation® software at a fixed distance of 55 cm. The visual angle was \textasciitilde{}1.07° for standard/distractor and \textasciitilde{}1.22° for target stimuli, shown on a dark gray background with light gray text. The Inter-Stimulus-Interval (ISI) was in the range from 1200 - 1600 ms and stimuli were separated by a blank screen. The total number of stimuli was 660, of which 12\% were targets, 12\% were deviants, and the remaining 76\% were standards. The participant was asked to press the key “2” on a numerical keyboard using the right middle finger in response to the target stimulus and refrain from responding to all other stimuli. The task was interspersed with four 15-second breaks during which a white fixation cross appeared in the center of a screen. The task lasted \textasciitilde{}22 minutes. The participant became familiar with the main task during a 3-minute training session consisting of \textasciitilde{}50 stimuli.

EEG data were recorded using the actiCHamp system and Brain Vision Recorder (Brain Products GmbH, Munich, Germany). A 128-channel actiCAP with active electrodes was used, with FCz as the online reference. Electrode positions were digitized using the CapTrak 3D scanner. Impedances were kept below \textasciitilde{}5 $k\Omega.$ Data were sampled at 1,000 Hz with a 280 Hz low-pass filter; no high-pass or notch filters were applied during recording. The EEG recording was conducted at the Nencki Institute of Experimental Biology, PAS, Warsaw, Poland, in a quiet, dimly lit room, with the participant seated comfortably facing a monitor. The researcher monitored the session remotely via desktop connection and LAN camera.

EEG data preprocessing was performed automatically using the MNE package \cite{MNE2014}. First, the signals were high-pass filtered at 1 Hz. Bad channels were identified based on abnormal power levels, defined as values exceeding 5 standard deviations (SD) in the 0-5 Hz range or 2.5 SD in the 5-40 Hz range. Transient artifacts with amplitudes greater than 500 $\mu$V were removed. The data were then low-pass filtered at 100 Hz and re-referenced to the average reference. Independent Component Analysis (ICA) was conducted using the extended infomax algorithm, with the number of components set to match the rank of the data. ICA components were automatically classified with ICLabel \cite{Li2022}, and any components not labeled as “brain” or “other” with a probability of at least 0.75 were discarded. Following ICA, the signals were low-pass filtered at 40 Hz, previously removed bad channels were interpolated, and the data were downsampled to 256 Hz. An automated artifact rejection procedure was applied, excluding data segments exceeding a threshold of 100 $\mu$V. Following visual inspection, this threshold was further reduced to 50 $\mu$V to remove residual transient artifacts. The data were then segmented into epochs $[-200,800]$ ms relative to the onset of the visual stimulus. Finally, DC detrending was applied, followed by baseline correction using the $[-200,0]$ ms pre-stimulus interval.

The 3-stimulus oddball paradigm is widely used for investigating ERP correlates of cognitive (context updating and working memory) processes, such as the P3b component, which is elicited by rare target stimuli requiring a behavioral response. The P3b component typically reaches peak amplitude at the Pz electrode between 350 and 600 ms after stimulus onset \cite{Polich2007} (see \ref{sensory}), and this interval has been selected as the ``cognitive'' interval.
In the visual 3-stimulus oddball paradigm, early sensory potentials could also be analyzed as the markers of perceptual and attentional processing that occur prior to the P3b component. The sensory (visual) potentials are typically observed over occipital electrodes (e.g., Oz, see \ref{sensory}) and include the N1 (N75) with the maximum amplitude \textasciitilde{}60-90 ms, and P1 (P100) component, which peaks \textasciitilde{}90-130 ms after the stimulus onset \cite{Creel2019}. Accordingly, we defined the 50-200 ms post-stimulus window as the ``sensory'' interval.

\subsubsection{Implementation Tutorial} \label{IT}
To conserve space, we show only the tutorial sections directly relevant to the \ref{localize!}-\ref{reconstruct!} pipeline introduced in Section~\ref{MNE_IMP}. The complete tutorial, including preprocessing steps and dataset link, is available at
\vspace{\baselineskip}
\\\url{https://julia-jurkowska.github.io/mvpure-tools}.
\vspace{\baselineskip}

\begin{lstlisting}[style=journalpython, caption={Empirical eigenvalue spectrum of $\tilde{\R}\tilde{\N}^{-1}$, where $\tilde{\R}$ and $\tilde{\N}$ denote finite sample estimates of $\R$ and $\N$, respectively.}]
fig_RN = viz.plot_RN_eigenvalues(
R=data_cov_task.data,
N=noise_cov.data,
subject=subject,
s=14
)
\end{lstlisting}

\begin{figure}[h!]
\centering
\includegraphics[width=0.9\linewidth]{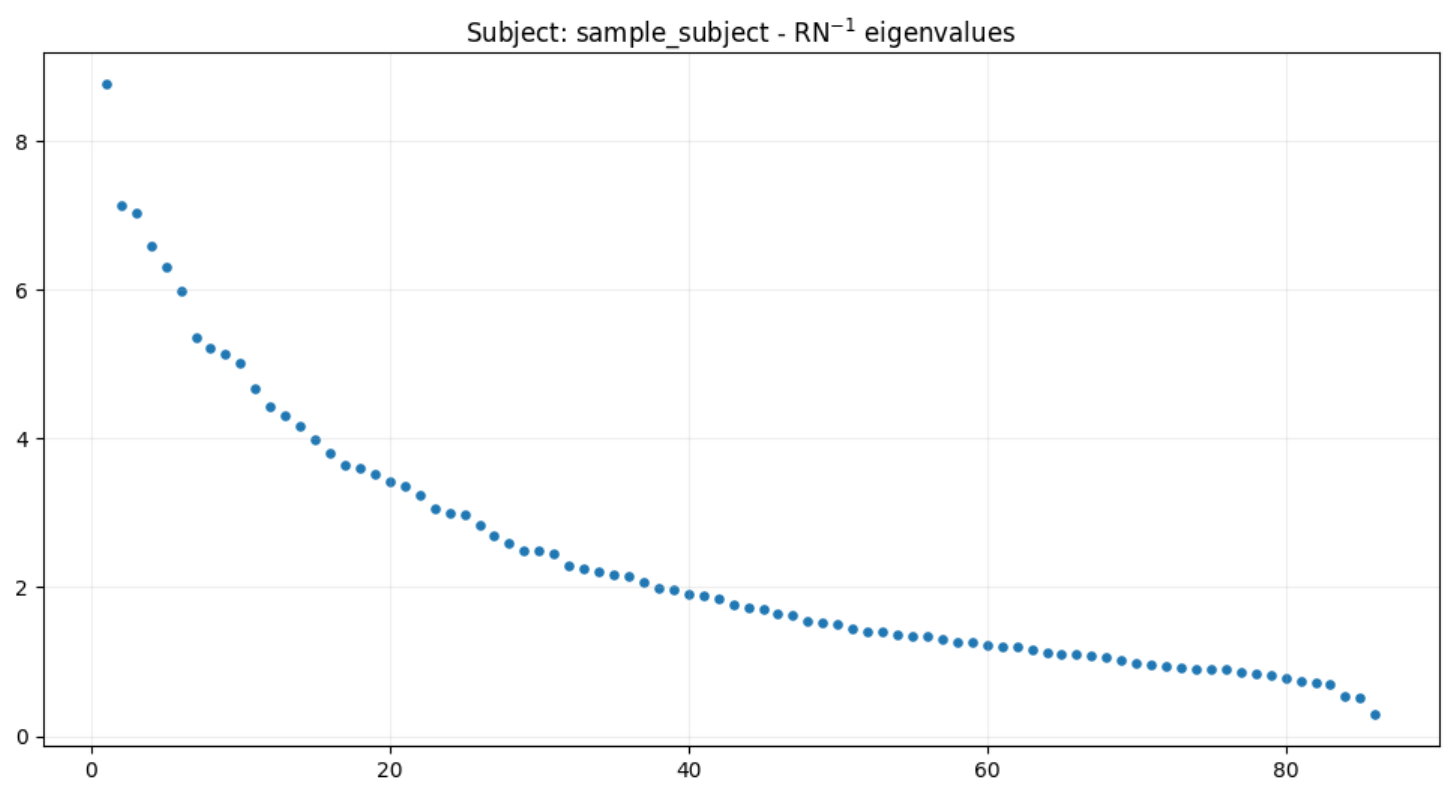}
\caption{Empirical distribution of eigenvalues of $\tilde{\R}\tilde{\N}^{-1}$ obtained from the finite sample estimates for the interval $[350,600]$ ms after target stimulus onset.} \label{eigs_fig}
\end{figure}

\noindent The above code generates Figure \ref{eigs_fig}, showing the empirical distribution of the eigenvalues of $\tilde{\R}\tilde{\N}^{-1}$, where $\tilde{\R}$ and $\tilde{\N}$ denote finite sample estimates of $\R$ and $\N$, respectively.

We note that the empirical distribution of the eigenvalues of $\tilde{\R}\tilde{\N}^{-1}$ deviates from the theoretical prediction in (\ref{eigs_RN}), where the number of eigenvalues strictly greater than one defines the dimensionality of the signal subspace. In practice, the empirical spectrum shows an increased spread of eigenvalues, a rather universal behaviour observed for general distributions \cite{Bai2010}. Further, the product structure $\mathbf{R}\mathbf{N}^{-1}$ further precludes a simple analytic characterization of the non-asymptotic eigenvalue behavior.

Nonetheless, the theoretical insights from (\ref{eigs_RN}) provide a principled basis for estimating the number of directions in the source space. Specifically, we identify gaps in the empirical eigenvalue spectrum of $\tilde{\R}\tilde{\N}^{-1}$, proceeding from the largest eigenvalues downward. These gaps delineate directions associated with higher source-level signal-to-noise ratios from less dominant components. Combined with evaluating $MAI_{MVP}$ across all ranks during localization, thereby fully exploiting its unbiasedness, this approach enables robust identification of the dominant drivers within the highest signal-to-noise ratio directions.

\begin{lstlisting}[style=journalpython, caption={Source localization implemented using the proposed multi-source neural activity indices (\ref{localize!}, Section \ref{MNE_IMP}).}]
locs_task = localizer.localize(
subject=subject,
subjects_dir=subjects_dir,
localizer_to_use=["mai_mvp"],
n_sources_to_localize=l_0,
R=data_cov_task.data,
N=noise_cov.data,
forward=fwd,
r=list(range(1, l_0 + 1))
)
\end{lstlisting}
In this example, we employ the $MAI_{MVP}$ neural activity index with $l_0 = 5$ directions. The localization procedure includes a flexible rank parameter $r$, specified either as a single value or as a set of candidate ranks (up to $l_0$). When a set is provided, separate localization runs are performed for each rank.

We follow the recommended strategy of evaluating all admissible ranks of the $MAI_{MVP}$ index. For each $r \in \{1, \ldots, l_0\}$, sources are localized using the corresponding rank-$r$ localizer and mapped onto cortical vertices. We apply the iterative localization algorithm described in \ref{localize!}, which produces at most $l_0(l_0+1)/2$ distinct candidate sources.

For the reconstruction \ref{reconstruct!}, the localized sources define a reduced forward model restricted to the selected vertices. Based on this subset, an LCMV beamformer is constructed using the MNE-Python implementation to estimate source time courses. No orientation constraint is imposed (\texttt{pick\_ori=None}), and Neural Activity Index (NAI) is applied to enhance the interpretability of the reconstructed activity.\footnote{Note that NAI is used here not as a localizer, but as a built-in normalization of the LCMV filter output.} The resulting spatial filter is then applied to the sensor-level data, yielding source-space time series at the localized vertices.

\begin{lstlisting}[style=journalpython, caption={Reconstruction of source activity for selected vertices using an LCMV beamformer with NAI normalization.}]
new_fwd = utils.subset_forward(
old_fwd=fwd,
localized=locs_task,
hemi='both'
)

filter_task = mne.beamformer.make_lcmv(
signal_task.info,
new_fwd,
data_cov_task,
reg=0.05,
noise_cov=noise_cov,
pick_ori=None,
weight_norm="nai",
rank=None
)

stc = mne.beamformer.apply_lcmv(signal_task, filter_task)
\end{lstlisting}

\subsubsection{Results}
\begin{figure}[H]
\centering
\begin{subfigure}[b]{\textwidth}
\centering
\includegraphics[width=0.65\linewidth]{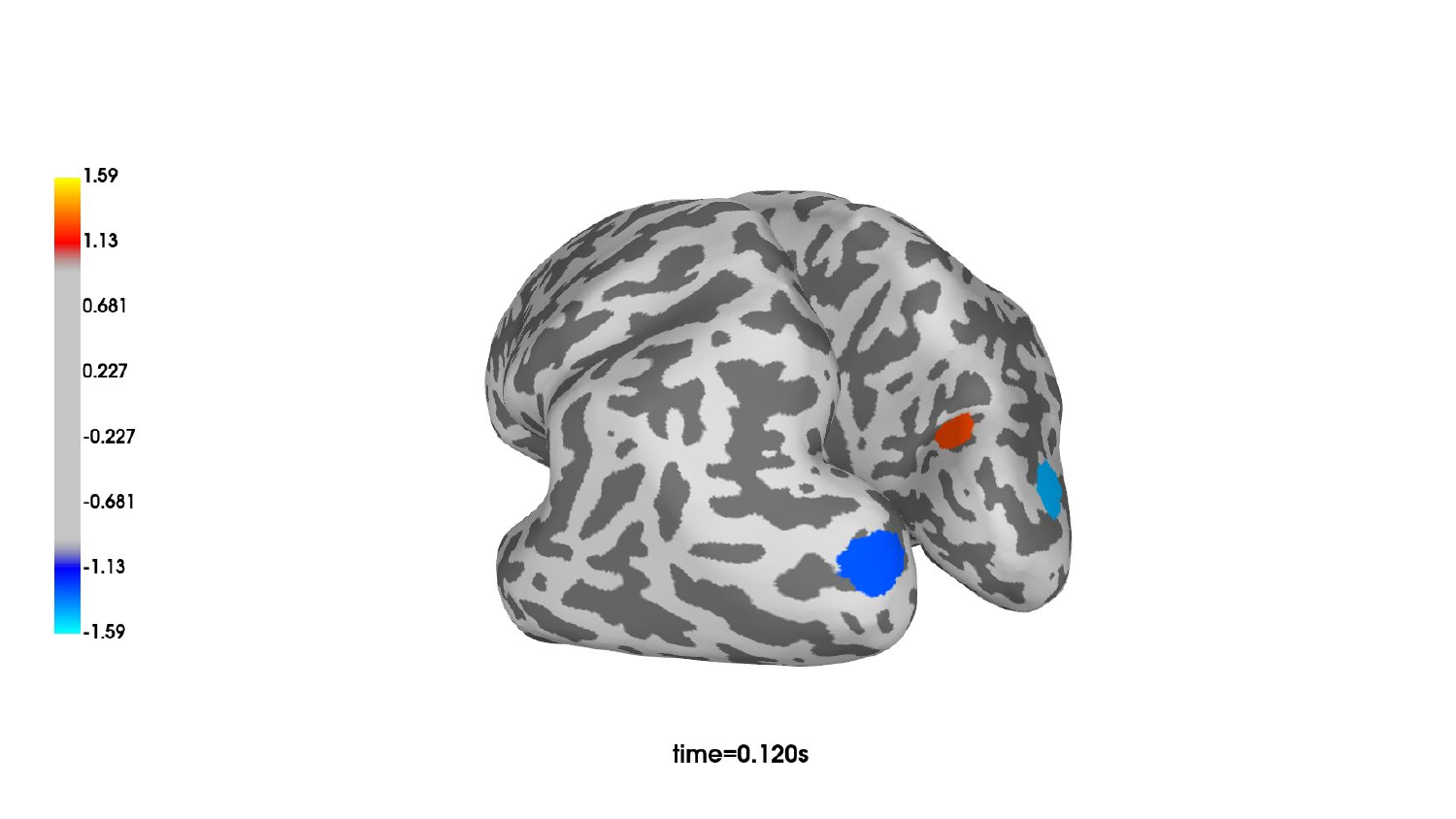}
\caption{A snapshot of the activity within the ``sensory'' interval.}
\label{fig:sub1}
\end{subfigure}
\vspace{0.5cm} 
\begin{subfigure}[b]{\textwidth}
\centering
\includegraphics[width=1.17\linewidth]{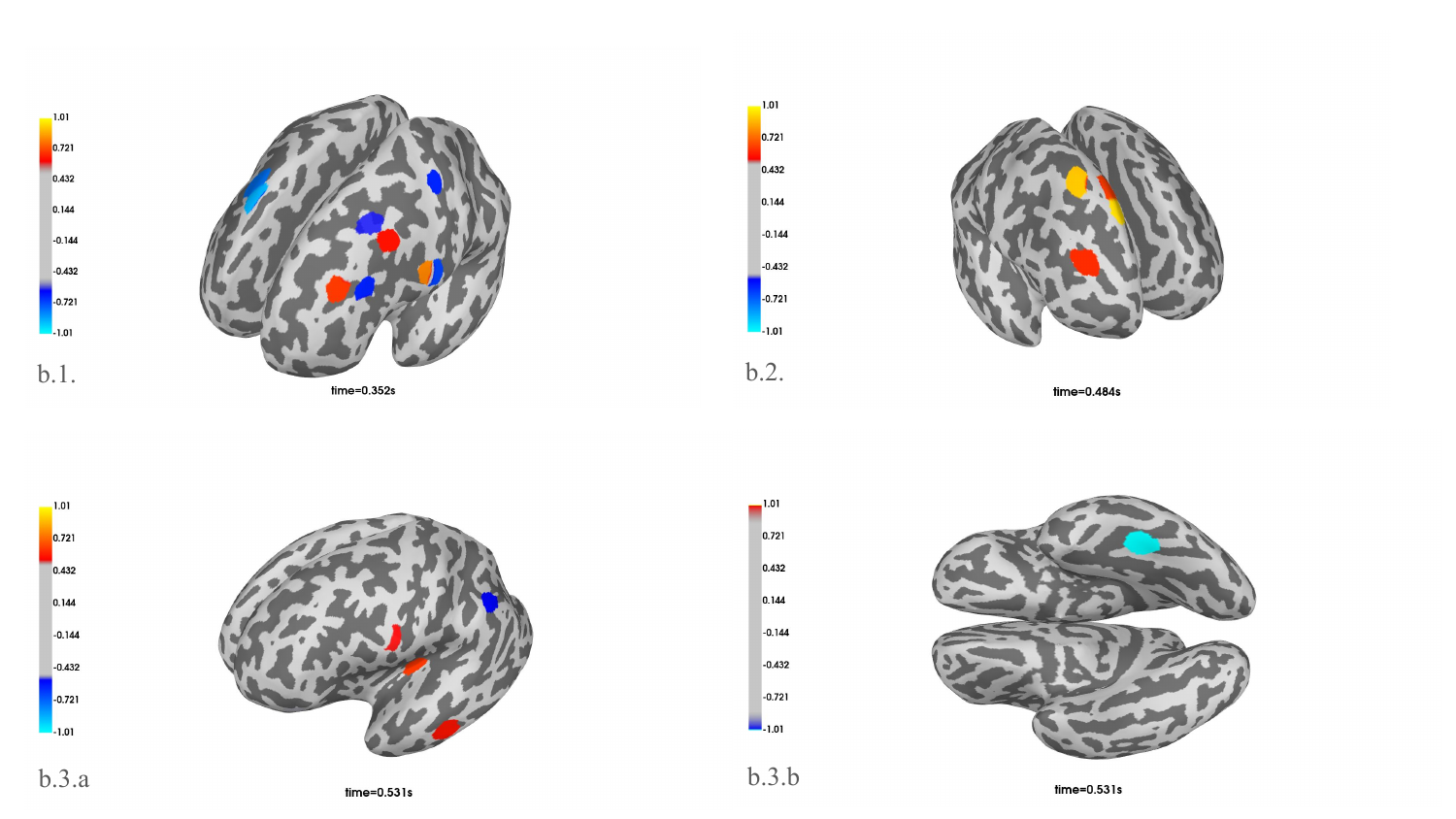}
\caption{A snapshot of the activity within the ``cognitive'' interval.}
\label{fig:sub2}
\end{subfigure}
\caption{EEG sources elicited by target stimuli in the visual oddball paradigm during early (``sensory'') and late (``cognitive'') time windows.
(a) Early sensory time window (50–200 ms; five-direction analysis), at 120 ms (panel a).
(b) Late cognitive time window (350–600 ms; ten-direction analysis), at 352 ms (b.1), 484 ms (b.2), and 531 ms (b.3.a, b.3.b). See \ref{sensory}, Tables \ref{S3}-\ref{C10}.}
\label{fig:two_subfigures}
\end{figure}

Consistent with existing literature on EEG source analysis in response to target stimuli presented in an active oddball paradigm \cite{Polich2007, Bledowski2004, Johnson1993, Linden2005,  Rusiniak2013, Sabeti2016, Tarkka1996}, a limited number of spatially localized sources were expected within the early ``sensory'' time window (50–200 ms post-stimulus onset), reflecting sensory processing and stimulus discrimination occurring locally within the brain (i.e., within primary and associative cortices). In contrast, multiple spatially proximal and temporally overlapping sources were expected within the 350–600 ms ``cognitive'' interval, associated with cognitive processing arising from interactions among distinct, albeit sometimes overlapping, neural networks (e.g., parietal-frontal regions and medial temporal areas). These time windows were defined based on our ERP results (\ref{sensory}, Figs.\ref{Oz} and \ref{Pz}).
We follow the recommended strategy of evaluating all admissible ranks of the $MAI_{MVP}$ index (see Section \ref{IT} for implementation tutorial). Source localization results were neurophysiologically plausible and consistent across ranks, indicating stability of the proposed pipeline. Source localization results are presented in Fig.\ref{fig:sub1}–\ref{fig:sub2}. A complete list of sources with corresponding MNI coordinates is provided in \ref{sensory}, Tables \ref{S3}-\ref{C10}.

\pagebreak

Within the ``sensory'' time window, both three- and five-direction analyses revealed the strongest sources at approximately 120 ms post-stimulus onset in the right superior occipital gyrus and the bilateral middle occipital gyri, corresponding to associative visual cortices (V2 and V3) (Tables \ref{S3} and \ref{S5}). The five-direction analysis further extends this set of sources by including the right cuneus, a region encompassing parts of the primary and secondary visual cortices (V1/V2) (Fig. 4, snapshot panel a; Table \ref{S5}). These findings are consistent with the expected involvement of associative visual cortical regions in early-stage stimulus processing \cite{Clark1994, Ungerleider1994}, thereby supporting the physiological plausibility of the localized sources.
Apart from activations in the visual cortices within the “sensory” interval, both the three- and five-direction analyses revealed the strongest sources in the left anterior cingulate and middle frontal gyri, as well as the right superior frontal gyrus (Tables \ref{S3} and \ref{S5}). The five-direction solution additionally revealed activations in the right middle cingulate gyrus and the left middle temporal and precentral gyri (Table \ref{S5}). The anterior and middle cingulate gyri, together with the bilateral frontal cortices, are associated with top-down modulation of visual processing, salience detection, and the recruitment of attentional mechanisms elicited by the stimulus \cite{Hickey2010, Menon2010, Paneri2017}. The left precentral gyrus may be involved in the preparation of motor responses to target stimuli (e.g., button pressing), whereas the left middle temporal gyrus may be engaged in early higher-order visual processing within the ventral stream and feature-based stimulus discrimination \cite{Ungerleider1994, Kravitz2013}. Overall, these results indicate that the proposed method effectively captures localized early-stage activity during an active visual oddball task while maintaining spatial specificity.
In contrast to the localized patterns observed within the “sensory” time window, both six- and ten-direction analyses performed within the “cognitive” interval, demonstrated a distributed network involving multiple cortical regions typically implicated in the detection of target stimuli in oddball tasks (e.g., \cite{Polich2007, Linden2005, Rusiniak2013, Linden1999}). These regions included predominantly frontal areas (bilateral superior and middle frontal cortices, left insula, as well as the precentral gyri and the left inferior frontal cortex), followed by temporo-occipital regions (bilateral fusiform and inferior temporal gyri, the right middle temporal gyrus, and the left parahippocampal gyrus), and parietal regions (the bilateral postcentral gyri, the left inferior and superior parietal lobules, the right angular gyrus, and the left precuneus) (Tables \ref{C6} and \ref{C10}). The strongest sources were observed at 352 ms (right superior frontal gyrus, left precentral gyrus, left inferior and middle frontal gyri, left insula, and left inferior temporal and parahippocampal gyri), 484 ms (right superior and medial frontal gyri), and 531 ms (left inferior frontal gyrus, left insula, left inferior parietal lobule, left fusiform gyrus, and left inferior temporal gyrus) post-stimulus onset (Figure 4, panel b.1, b.2, b.3a, b.3b).
In the present study, the involvement of multiple superior, middle, and inferior frontal regions during the cognitive time window is consistent with their role in higher-order cognitive control processes underlying target detection and response selection \cite{Barbey2013}. The left insula may contribute to pre-attentive detection of deviant stimuli, signaling the need for attentional reorientation and engagement of higher-order cognitive networks \cite{Menon2010}. Additionally, in this time window the insula may be involved in integration of interoceptive and sensory information to guide adaptive responses to novel stimuli \cite{Uddin2015}.
The precentral and postcentral gyri likely reflect motor preparation and execution associated with the behavioral response.
The sources identified in the temporal cortices (fusiform, inferior temporal and parahippocampal gyri) form part of the ventral visual stream, which, together with medial temporal structures, supports higher-order perceptual and mnemonic processing. In the active oddball task, these regions may contribute to detailed perceptual analysis and identification of target stimuli, integration of visual features into coherent object representations, and updating of contextual representations of these stimuli in memory \cite{Kravitz2013, Grill2004, Epstein1998}.
Parietal regions, including the bilateral inferior and superior parietal lobules, as well as the left precuneus and the right angular gyrus, may reflect goal-directed attentional modulation, context updating, and evaluation of stimulus significance, all of which are required for appropriate evaluation of and response to target stimuli \cite{Cavanna2006, Corbetta2002}.
Overall, the results demonstrate that the proposed method is capable of resolving multiple spatially and temporally overlapping sources in real EEG data while maintaining neurophysiological plausibility. Notably, similar spatial patterns were observed across different dimensions (e.g., three- versus five-direction analyses in the “sensory” window and six- versus ten-direction analyses in the “cognitive” window), indicating the robustness of the localization with respect to the assumed source dimensionality.

\section{Discussion} \label{DISCUSSION}
We proposed a two-stage localization-reconstruction framework for the EEG/MEG inverse problem, designed for scenarios in which multiple focal sources dominate the measured data. The framework is operationalized using novel multi-source neural activity indices in \ref{localize!}, coupled with an MNE-Python-based implementation of the LCMV spatial filter in \ref{reconstruct!}. A representative example is task-evoked activity, such as in oddball paradigms, where multiple focal, functionally coupled and potentially correlated cortical sources are simultaneously active and jointly contribute to the measured signals. The use-case tutorial in Section \ref{oddball} demonstrates how the proposed method can be implemented by extending MNE-Python in a way that is natural for its users, with the goal of facilitating adoption within the neuroscience community. We emphasize that this use case focuses on a surface lead-field model and EEG data, settings that have been comparatively underexplored for beamformer-based methods. Furthermore, thanks to its MNE-Python-based implementation, the method is readily extensible. In the following, we outline natural extensions that merit consideration for future research.

First, the reconstruction stage \ref{reconstra} can be extended beyond the MNE-Python-based LCMV filter used in \ref{reconstruct!} to a broader class of multi-source filters, including, but not limited to, the reduced-rank multi-source MV-PURE filters introduced in \cite{Piotrowski2019}. Recent comparative studies have demonstrated that MV-PURE-based filters achieve strong performance in multi-source scenarios, including when combined with other multi-source-aware approaches \cite{Hecker2026}, with benchmark results available at
\\\url{https://lukethehecker.github.io/invertmeeg/benchmarks/}.

\noindent This line of work is both promising and significant, particularly in the context of integrating such filters within the MNE-Python ecosystem, where they can leverage its infrastructure and capabilities, as does the existing LCMV implementation.

Second, a dedicated group-level visualization tool is essential for any study investigating source-level activity at the group level. Simply overlaying single-subject activations onto an average brain is inadequate: inter-subject anatomical variability and individual patterns of activation mean that spatial overlays can obscure rather than reveal common effects. Instead, methods for robustly selecting voxels and regions that contribute the most relevant activity across subjects - together with transparent ranking schemes that quantify each contributor’s relative importance - are needed. Such tools should combine aggregation (e.g., effect-size or weighted statistical maps), anatomically informed parcellation, and visualization of uncertainty and variability to support group inferences.

Third, extending the implementation of the neural activity indices to volumetric, orientation-unconstrained (vector) forward models is necessary and especially valuable for the MEG inverse problem. Mathematically, the current repository already supports vector lead-fields: each candidate source is expanded by a unit orientation vector, and the index is evaluated over the orientation(s) that maximize the neural activity measure. However, this substantially enlarges the search space, as its dimensionality becomes $s\times p$ (with $s$ candidate sources of $p$ candidate orientations). Fortunately, evaluation of the proposed neural activity indices is naturally parallelizable and, in principle, limited only by available computational resources, since calls to the neural activity index evaluation function are independent and can be executed in parallel. Accordingly, CPU- and GPU-based parallel implementations of the indices in \ref{localize!} are feasible and should be pursued as a future research.

Fourth, a comprehensive performance evaluation of the proposed two-stage localization-reconstruction pipeline is needed, given that the current results position the method alongside widely adopted M/EEG inverse solutions in terms of completeness and practical utility. While previous simulation studies have investigated both neural activity indices, e.g., the full-rank MAI \cite{Moiseev2011} and the ``first-generation'' reduced-rank indices \cite{Piotrowski2021}, as well as MV-PURE-based multi-source spatial filters \cite{Piotrowski2019}, the present stage of methodological development calls for a systematic comparison with state-of-the-art solutions commonly used in the neuroscientific community. Of particular interest is whether the proposed approach, leveraging the full spectrum of ranks of the multi-source neural activity indices (\ref{localize!}), provides increased robustness to forward model inaccuracies, especially when coupled with the multi-source spatial filters to be implemented in the reconstruction stage (\ref{reconstruct!}).

Finally, recent developments in machine learning for EEG/MEG inverse problems suggest it may become computationally feasible to relax the discrete rank constraint used in the neural activity indices in \ref{localize!} by replacing it with continuous regularization using the convex nuclear norm. However, whether the additional computational cost of this relaxation yields meaningful improvements in spatial resolution remains an open question.

\section{Ethics statement} \label{ethics}
The data analyzed in Section \ref{exp_data} were obtained as part of a previously approved research project. All participants in the original study \cite{Dzianok2022} provided written informed consent for participation, including consent for the anonymized use of their data in subsequent analyses and scientific publications. This study was conducted in accordance with the Declaration of Helsinki and was covered by the approval of the Research Ethics Committee, Faculty of Humanities\footnote{At present, the Ethics Committee is affiliated with the Faculty of Philosophy and Social Sciences, Nicolaus Copernicus University in Toruń, Poland (Gagarina 39, 87-100 Toruń, Poland). \url{https://www.wfins.umk.pl/wydzial/komisja-ds-etyki-badan-naukowych/}}

For the purposes of the present study, a single case was selected from the existing dataset. All data were fully anonymized prior to analysis, and no information allowing identification of the participant is disclosed. Participants in the original study received financial compensation (200 PLN, ~44 EUR).

\section*{Acknowledgments}
The authors would like to thank the Excellence Initiative - Research University programme at Nicolaus Copernicus University in Toruń: the University Centre of Excellence ``Dynamics, Mathematical Analysis and Artificial Intelligence'' and the Emerging Field ``Culture, Development and Wellbeing''. The authors would also like to thank Mr. Jan Nikadon and Dr. Michał Lemańczyk for helpful discussions. Finally, the authors are grateful to anonymous reviewers for their constructive comments which surely promoted the readability of the revised manuscript.

\pagebreak
\appendix

\section{Iterative Source Search Algorithm} \label{ISSA}
A brute-force strategy to evaluate a multi-source neural activity index for $l$ sources requires $\binom{s}{l}$ evaluations (with $s$ the number of candidate sources) and is computationally infeasible for realistic $s$. To overcome this combinatorial bottleneck we exploit two features of the proposed approach. First, Definition~\ref{DFMVP} introduces the $MAI_{MVP}$ indices in forms that are defined for any candidate source count~$l$ and that naturally admit an incremental search: one may discover sources sequentially by starting with $l=1$ and then adding sources until $l=l_0.$ Second, Proposition~\ref{FULLFORM2025} provides compact, computationally efficient representations of these indices that make the incremental updates inexpensive to evaluate. Together, these results replace an intractable combinatorial search by a practical iterative localisation procedure: at each step the algorithm evaluates the index for only $s$ candidate additions and selects the maximiser, repeating until $l_0$ sources have been found. Consequently, the search cost is reduced from combinatorial $\binom{s}{l_0}$ evaluations to a sequence of evaluations that scales roughly linearly with $s$ per discovered source.

The iterative source-selection scheme employed here builds directly on earlier work \cite{Moiseev2011} which demonstrated its efficient implementation for the multi-source indices, and \cite{Piotrowski2021} which  presented its heuristic adaptation for the ``first-generation'' reduced-rank indices whose forms depend on \(\mathbf{Q}_0\) - the covariance matrix of \(\mathbf{q}_0\) in model~(\ref{model}), which is typically unknown. In contrast, Algorithm~\ref{A1} illustrates the iterative source-search applied to the unbiased reduced-rank multi-source neural activity indices \(MAI_{MVP}\), whose forms do not depend on $\Q_0.$

\begin{algorithm}[h] 
  \caption{Iterative discovery of sources using neural activity indices}
  \label{A1}
  \begin{algorithmic}[1]
    \REQUIRE $l$ - the number of sources considered, $r$ - the rank of the index, $s$ - the number of candidate sources, $MAI_{MVP}$ neural activity indices
    \\ \textit{Initialisation}
    \STATE $\bth=\{\varnothing\}$
    \\ \textit{LOOP Process}
    \FOR {$j=1$ to $l$}
    \FOR {$i=1$ to $s$}
    \STATE set $\bth_i=\bth\cup\{\theta_i\}$
    \ENDFOR
    \IF {$j\leq r$}
    \STATE $\bth=\underset{i=1,\dots,s}{\arg\max}\ \sum_{i=1}^{l}\lambda_i(\G(\bth_i)\S(\bth_i)^{-1})-l$ [see (\ref{MAIMVP_small})]
    \ELSE
    \STATE $\bth=\underset{i=1,\dots,s}{\arg\max}\ \sum_{i=1}^{r}\lambda_i(\G(\bth_i)\S(\bth_i)^{-1})-r$ [see (\ref{MAIMVP}) and (\ref{MAIMVPEXT})]
    \ENDIF
    \ENDFOR
    \RETURN $\bth$
  \end{algorithmic}
\end{algorithm}

\pagebreak
\section{Known Results Used} \label{kru}
\begin{fact} \label{Q_EST}
Under assumptions of model \ref{model} and with notation as in (\ref{QN}), (\ref{R}), (\ref{G_Filt}) and (\ref{S_Filt}), one has
\begin{equation}
\Q_0=\S_0^{-1}-\G_0^{-1}.
\end{equation}
The proof is provided in, e.g., \cite[Lemma~1]{Piotrowski2014}.
\end{fact}

\begin{fact}[\cite{Theobald1975}] \label{Theobald}
Let $\C\in\sgen{n}{n},\D\in\sgen{n}{n}$ be symmetric matrices. Denoting by $c_1\geq c_2\geq\dots\geq c_n$ and $d_1\geq d_2\geq\dots\geq d_n$ the eigenvalues of $\C$ and $\D$, respectively, one has 
\begin{equation} \label{poeq}
tr\{\C\D\}\leq\sum_{i=1}^n c_id_i.
\end{equation}
\end{fact}

\begin{fact}[{\citep[Theorem 7.7.2]{Horn1985}}] \label{Horn772}
Let $\A,\B\in\sgen{n}{n}$ be symmetric and let $\Ss\in\sgen{n}{m}.$ One has
\begin{equation*}
\A\succeq\B\implies\Ss^t\A\Ss\succeq\Ss^t\B\Ss,
\end{equation*}
and
\begin{equation*}
rank(\Ss)=m\textrm{ and }\A\succ\B\implies\Ss^t\A\Ss\succ\Ss^t\B\Ss.
\end{equation*}
\end{fact}
  
\begin{fact}[{\citep[Corollary 7.7.4]{Horn1985}}] \label{Horn774}
Let $\A,\B\in\sgen{n}{n}$ be symmetric. Let $\la_i(\A),\ \la_i(\B)$ be the ordered eigenvalues of $\A$ and $\B$, respectively. One has
\begin{equation*}
(\A\succ\0\textrm{ and }\B\succ\0)\implies(\A\succeq\B\iff\B^{-1}\succeq\A^{-1}),
\end{equation*}
and
\begin{equation*}
\A\succeq\B\implies\la_i(\A)\geq\la_i(\B),\ i=1,\dots,n.
\end{equation*}
\end{fact}

\begin{fact}[{\citep[Observation 7.1.8]{Horn1985}}] \label{Horn718}
Let $\A\in\sgen{n}{n}$ be symmetric and let $\C\in\sgen{n}{m}.$ Then $rank(\C^t\A\C)=rank(\A\C).$
\end{fact}

\section{Proof of Proposition \ref{FULLFORM2025}} \label{PD_FULLFORM2025}
\ref{MVP_unbiased}. 
The proof of \ref{MVP_unbiased} adapts parts of the proofs of \cite[Lemma 1]{Piotrowski2021} and \cite[Theorem 1]{Piotrowski2021} amended to the current setup below. We recall that matrices which depend on the true sources' locations $\bth_0$ are denoted as $\H_0=\H(\bth_0)$, $\G_0=\G(\bth_0)$, and so on. Similarly, to simplify derivation below, we omit the source location parameter vector $\bth$ wherever it does not lead to ambiguity, so $\H(\bth)$ is replaced with $\H$, $\G(\bth)$ is replaced with $\G$, and so on.

Let $\I^r_s\in\sgen{s}{s}$ be a matrix with its $r\times r$ principal submatrix the identity matrix and zeros elsewhere if $r<s$ and the $\I_s$ identity matrix if $r\geq s$, and let 
\begin{equation} \label{tildeH}
\tilde{\H}\df \H\Q^{1/2}\in\sgen{m}{l}.
\end{equation}
Further, let $\U\Si\V^{t}$ be the singular value decomposition of $\R^{-1/2}\tilde{\H}$ with diagonal entries of $\Si$ organized in non-increasing order. Then $\tilde{\H}=\R^{1/2}\U\Si\V^t$ and hence, for $l>r$,
\begin{multline} \label{bound} 
MAI_{MVP}(\bth,r)+r=tr\{\V\Si^{t}\U^{t}\R^{1/2}\N^{-1}\R^{1/2}\U\Si\V^{t}\\
(\V\Si^{t}\U^{t}\U\Si\V^{t})^{-1}\V\I^{r}_{l}\V^{t}\}=\\
tr\{\Si^{t}\U^{t}\R^{1/2}\N^{-1}\R^{1/2}\U\Si(\Si^{t}\Si)^{-1}\I^{r}_{l}\}=\\
tr\{\R^{1/2}\N^{-1}\R^{1/2}\U\Si(\Si^{t}\Si)^{-1}\I^{r}_{l}\Si^{t}\U^{t}\}=\\
tr\{\R^{1/2}\N^{-1}\R^{1/2}\U\I^{r}_{m}\U^{t}\}
\leq\sum\limits_{i=1}^{r}\lambda_{i}(\R^{1/2}\N^{-1}\R^{1/2})=\\
\sum\limits_{i=1}^{r}\lambda_{i}(\R\N^{-1}).
\end{multline}
The inequality above is obtained from Theobald's theorem (Fact \ref{Theobald} in \ref{kru}) and the last equality from matrix similarity and the fact that  $\R^{1/2}\N^{-1}\R^{1/2}\succ\0$ by Fact \ref{Horn772} in \ref{kru}. For $l\leq r$, one can reproduce the above derivation to arrive at
\begin{equation} \label{obv1}
MAI_{MVP}(\bth,r)+l\leq\sum\limits_{i=1}^{l}\lambda_{i}(\R\N^{-1})\leq\sum\limits_{i=1}^{r}\lambda_{i}(\R\N^{-1}).
\end{equation}
Therefore, since we are interested in maximizers of the introduced indices, we restrict the subsequent derivation of unbiasedness of $MAI_{MVP}$ to the case $l>r.$

Now, to prove unbiasedness of $MAI_{MVP}$, we use Fact \ref{Q_EST} in \ref{kru}, so that
\begin{equation} \label{s0g0}
\tilde{\S}_0^{-1}=\Q_0^{-1/2}\S_0^{-1}\Q_0^{-1/2}=\Q_0^{-1/2}\G_0^{-1}\Q_0^{-1/2}+\I_{l_0}=\tilde{\G}_0^{-1}+\I_{l_0},
\end{equation}
where $\tilde{\G}_0^{-1}=\Q_0^{-1/2}\G_0^{-1}\Q_0^{-1/2}.$ Thus, in particular, 
\begin{equation} \label{proc}
\Pre^{(r)}_{\tilde{\S}_0}=\Pre^{(r)}_{\tilde{\G}_0}.
\end{equation}
Therefore, by (\ref{s0g0}) and (\ref{proc})
\begin{multline} \label{unbiased!} MAI_{MVP}(\bth_{0},r)=tr\{\tilde{\G}_0(\tilde{\S_0})^{-1}\Pre^{(r)}_{\tilde{\S_0}}\}-r=\\
  tr\{\tilde{\G}_0(\tilde{\G}_0^{-1}+\I_{l_0})\Pre^{(r)}_{\tilde{\G_0}}\}-r=
  tr\{\tilde{\G}_0\Pre^{(r)}_{\tilde{\G_0}}\}=\sum_{i=1}^r\la_i(\tilde{\G}_0).
\end{multline}
Let $\lambda>0$ be an eigenvalue of $\tilde{\G}_0\succ\0$, i.e., there exists $0\neq\x\in\sgens{l_0}$ such that
$$\tilde{\G}_0\x=\Q_0^{1/2}\G_0\Q_0^{1/2}\x=\tilde{\H}_0^t\N^{-1}\tilde{\H}_0\x=\la\x,$$ thus also 
$$\tilde{\H}_0\tilde{\H}_0^t\N^{-1}\tilde{\H}_0\x=\la\tilde{\H}_0\x,$$
hence $\la>0$ is an eigenvalue of $\tilde{\H}_0\tilde{\H}_0^t\N^{-1}=\H_0\Q_0\H_0^t\N^{-1}.$ Since $\N\succ\0$ is invertible, we have that
\begin{equation} \label{rank1}
rank(\H_0\Q_0\H_0^t\N^{-1})=rank(\H_0\Q_0\H_0^t),
\end{equation}
and, from Fact \ref{Horn718} in \ref{kru}, we have that
\begin{equation*}
rank(\H_0\Q_0\H_0^t)=rank(\Q_0\H_0^t).
\end{equation*}
Moreover, since $\Q_0\succ\0$ is invertible, we conclude from the above that 
\begin{equation*}
rank(\H_0\Q_0\H_0^t\N^{-1})=rank(\H_0^t)=l_0.
\end{equation*}
Further, $\H_0\Q_0\H_0^t\N^{-1}$ is similar to $\N^{-1/2}\H_0\Q_0\H_0^t\N^{-1/2}\succeq\0$ by Fact \ref{Horn772} in \ref{kru}, hence it has $l_0$ positive eigenvalues and $m-l_0$ zero eigenvalues. We conclude from the above that, if $\lambda>0$ is an eigenvalue of $\tilde{\G}_0$, then it is also one of the $l_0$ positive eigenvalues of $\H_0\Q_0\H_0^t\N^{-1}.$ We recall now that $$\R=\H_0\Q_0\H_0^t+\N,$$ hence
\begin{equation} \label{rank_fin}
\H_0\Q_0\H_0^t\N^{-1}=\R\N^{-1}-\I_m.
\end{equation}
Therefore, from (\ref{unbiased!}) and the above considerations we conclude that
\begin{equation} \label{unbiased!!}
MAI_{MVP}(\bth_{0},r)=\sum_{i=1}^r\la_i(\tilde{\G}_0)=\sum_{i=1}^r\la_i(\R\N^{-1}-\I_{l_0})=\sum_{i=1}^r\la_i(\R\N^{-1})-r.
\end{equation}
It now follows from (\ref{bound})-(\ref{obv1}) and (\ref{unbiased!!}) that $MAI_{MVP}(\bth_{0},r)$ is unbiased (in the sense of Definition \ref{nai}) for any given $1\leq r\leq l_0.$

\ref{MVP_MAX}. We note first that, for $l>r$, we have:
\begin{multline} \label{KEY2025}
MAI_{MVP}(\bth,r)+r=tr\{\tilde{\G}\tilde{\S}^{-1}\Pre^{(r)}_{\tilde{\S}}\}=\\
tr\{\G\S^{-1}\Q^{-1/2}\Pre^{(r)}_{\tilde{\S}}\Q^{1/2}\}=
tr\{\G\S^{-1}\Pre^{(r)}_{\S\Q}\},
\end{multline}  
where
\begin{equation} \label{P_SQ}
\Pre^{(r)}_{\S\Q}\df \Q^{-1/2}\Pre^{(r)}_{\tilde{\S}}\Q^{1/2}
\end{equation}
is an oblique projection matrix. Now, from (\ref{tQ}) we obtain
\begin{equation} \label{SQGS}
\Q^{-1/2}\tilde{\S}\Q^{1/2}=\S\Q=\S(\S^{-1}-\G^{-1})=\I_l-\S\G^{-1}.
\end{equation}
Let $\V\Ga\V^t$ be an eigenvalue decomposition of $\tilde{\S}\succ\0$ with eigenvalues organized in nonincreasing order, so that
\begin{equation} \label{Prs}
\Pre^{(r)}_{\tilde{\S}}=\V\I^r_l\V^t,
\end{equation}
and, from (\ref{SQGS}),
\begin{equation} \label{SQ}
\S\Q=\Q^{-1/2}\V\Ga\V^t\Q^{1/2}
\end{equation}
is an eigenvalue decomposition of $\S\Q.$ Let further
\begin{equation} \label{GS}
\G\S^{-1}=\U\Si\U^{-1}
\end{equation}
be en eigenvalue decomposition of $\G\S^{-1}$, with positive eigenvalues organized in nonincreasing order: the existence of such a decomposition of $\G\S^{-1}$ is implied by the matrix similarity of $\G\S^{-1}$ to the semi-positive matrix $\S^{-1/2}\G\S^{-1/2}\succ\0.$ From (\ref{SQGS}), (\ref{SQ}), and (\ref{GS}), we obtain therefore that
\begin{equation} \label{KEY_MAR_2026}
\Q^{-1/2}\V\Ga\V^t\Q^{1/2}=\S\Q=\I_l-\S\G^{-1}=\U(\I_l-\Si^{-1})\U^{-1},
\end{equation}
thus, we may set $\U=\Q^{-1/2}\V$ and $\Si=(\I_l-\Ga)^{-1}$ in (\ref{GS}).\footnote{We learn en passant that all eigenvalues of both $\tilde{\S}$ and $\S\Q$ are in $(0,1).$} From the above considerations we conclude that
\begin{multline*}
MAI_{MVP}(\bth,r)+r=tr\{\G\S^{-1}\Pre^{(r)}_{\S\Q}\}=\\
tr\{\Q^{-1/2}\V(\I_l-\Ga)^{-1}\V^t\Q^{1/2}\Q^{-1/2}\V\I^r_l\V^t\Q^{1/2}\}=\\
tr\{(\I_l-\Ga)^{-1}\I^r_l\}=\sum_{i=1}^r(1-\ga_i)^{-1}=
\sum_{i=1}^r\la_i(\G\S^{-1}),
\end{multline*}
which completes the proof of (\ref{MAIMVPEXT}). $\blacksquare$

\pagebreak
\section{Supplementary Simulation Data} \label{sim_sup}
\begin{table}[h!]
\centering
\footnotesize
\setlength{\tabcolsep}{3pt}
\begin{tabular}{|p{4cm}|p{10cm}|}
\hline
\textbf{Parameter name} & \textbf{List of labels used} \\
\hline
all\_labels\newline (Both Simulations) &
\begin{tabular}[t]{@{}l@{}}
cuneus-lh\\
cuneus-rh\\
lateraloccipital-lh\\
lateraloccipital-rh\\
inferiorparietal-lh\\
inferiorparietal-rh\\
superiorparietal-lh\\
superiorparietal-rh\\
superiortemporal-lh\\
superiortemporal-rh\\
supramarginal-lh\\
supramarginal-rh\\
superiorfrontal-lh\\
superiorfrontal-rh\\
insula-lh\\
insula-rh\\
caudalanteriorcingulate-lh\\
caudalanteriorcingulate-rh
\end{tabular} \\
\hline
poststimuli\_labels\newline (Simulation 1) &
\begin{tabular}[t]{@{}l@{}}
lateraloccipital-lh\\
lateraloccipital-rh
\end{tabular} \\
\hline
poststimuli\_labels\newline (Simulation 2) &
\begin{tabular}[t]{@{}l@{}}
cuneus-lh\\
cuneus-rh\\
lateraloccipital-lh\\
lateraloccipital-rh\\
inferiorparietal-lh\\
inferiorparietal-rh
\end{tabular} \\
\hline
\end{tabular}
\caption{Anatomical labels used as active and background regions in each simulation.}
\label{tab:used_labels_desc}
\end{table}

\pagebreak
\begin{table}[h!]
\centering
\footnotesize
\setlength{\tabcolsep}{3pt}
\captionsetup{justification=raggedright,singlelinecheck=false}
\begin{tabular}{|p{5cm}|p{6.5cm}|p{3cm}|}
\hline
\textbf{Parameter name} & \textbf{Parameter description} & \textbf{Set value} \\
\hline
n\_epochs & Number of epochs to simulate in single simulation & 100 \\
\hline
n\_vertices\_per\_label\_bg & Number of vertices to randomly select within background activity labels & 2 (Scenario 1), \newline 1 (Scenario 2) \\
\hline
erp\_factor & Scaling factor for ERP waveform amplitude & 15e-9 \\
\hline
order\_bg & MVAR order for background activity & 7 \\
\hline
order & MVAR order for post stimuli activity & 3 \\
\hline
target\_std\_bg & Target standard deviation of background source activity in Am & 15e-9 \\
\hline
target\_std & Target standard deviation of post stimuli activity in Am & 15e-9 \\
\hline
coupling\_bg & Cross-source coupling strength for background MVAR model & 0.5 \\
\hline
coupling & Cross-source coupling strength for stimulus-driven MVAR model & 0.8 \\
\hline
noise\_bg & Innovation noise amplitude for background activity & 1.0 \\
\hline
noise & Innovation noise amplitude for stimulus-driven activity & 1.0 \\
\hline
gaussian\_filter\_sigma\_bg & Temporal smoothing for background activity & 3.0 \\
\hline
gaussian\_filter\_sigma & Temporal smoothing for post stimulus activity & 3.0 \\
\hline
n\_dominant\_eigvals & Dimensionality of dominant latent sources space used to generate correlated stimulus-locked activity & 2 \\
\hline
noise\_factor & Scaling factor controlling sensor noise amplitude & 0.3 \\
\hline
\end{tabular}
\caption{Parameters used during simulation experiments.}
\label{tab:constant-parameters-simulation}
\end{table}

\pagebreak
\section{Supplementary Experimental Data} \label{sensory}
We provide links to source reconstruction movies for both the ``sensory'' $[50,200]$ ms and the ``cognitive'' $[350,600]$ ms intervals:
\vspace{\baselineskip}\\\url{https://figshare.com/projects/EEG_MEG_Multi-Source_Spatial_Filters_with_mvpure-py_Data_Movies_/263188}.
\vspace{\baselineskip}

\noindent We also note that a complete, interactive analysis using MNE-Python source browsing functionality can be obtained by running the provided tutorial code, available at
\vspace{\baselineskip}
\\\url{https://julia-jurkowska.github.io/mvpure-tools}.
\vspace{\baselineskip}

Localization results (\ref{localize!}) for Section \ref{exp_data} are summarized in Tables X–Y. Coordinates are reported in Montreal Neurological Institute (MNI) space. Approximate anatomical labels were assigned using standard neuroanatomical atlases, with Brodmann areas (BA) indicated where applicable. Sources are ranked in descending order according to NAI-normalized dipole moment amplitudes obtained during reconstruction \ref{reconstruct!}.

\begin{table}[h!]
\centering
\footnotesize
\begin{tabular}{|c|c|p{0.55\linewidth}|}
\hline
\textbf{Source No.} & \textbf{MNI coordinates} & \textbf{Approx. Anatomical Location} \\ \hline
1 & 26.16, -90.30, 20.70 & Right superior occipital gyrus \\ \hline
2 & -18.12, -93.28, 9.67 & Left middle occipital gyrus \\ \hline
3 & 43.06, -83.34, 3.61 & Right middle occipital gyrus \\ \hline
4 & -5.17, 15.71, 30.82 & Left anterior cingulate gyrus (BA 24) \\ \hline
5 & 24.42, 31.74, 48.04 & Right superior frontal gyrus (BA 8) \\ \hline
6 & -25.95, 9.08, 48.59 & Left middle frontal gyrus \\ \hline
\end{tabular}
\caption{Sources in the ``sensory'' time window (50-200 ms) for target stimuli, ranked by normalized dipole moment amplitude (three-direction analysis, active visual oddball task, one participant). Brodmann areas are indicated in parentheses, where applicable.}
\label{S3}
\end{table}
\pagebreak

\begin{table}[h!]
\centering
\footnotesize
\begin{tabular}{|c|c|p{0.55\linewidth}|}
\hline
\textbf{Source No.} & \textbf{MNI coordinates} & \textbf{Approx. Anatomical Location} \\ \hline
1 & 26.16, -90.30, 20.70 & Right superior occipital gyrus \\ \hline
2 & 15.86, -75.75, 22.86 & Right cuneus (BA 18) \\ \hline
3 & -18.12, -93.28, 9.67 & Left middle occipital gyrus \\ \hline
4 & 25.20, -75.61, 13.60 & Right cuneus \\ \hline
5 & 43.06, -83.34, 3.61 & Right middle occipital gyrus \\ \hline
6 & -5.17, 15.71, 30.82 & Left anterior cingulate gyrus (BA 24) \\ \hline
7 & 6.97, -33.62, 29.81 & Right middle cingulate gyrus \\ \hline
8 & 12.73, 35.94, 42.35 & Right superior medial frontal gyrus \\ \hline
9 & 24.42, 31.74, 48.04 & Right superior frontal gyrus (BA 8) \\ \hline
10 & -54.44, -2.89, -27.90 & Left middle temporal gyrus (BA 20) \\ \hline
11 & -25.95, 9.08, 48.59 & Left middle frontal gyrus \\ \hline
12 & -23.26, 64.98, -4.29 & Left superior frontal gyrus, orbital part (BA 10) \\ \hline
13 & -34.09, -4.15, 39.04 & Left precentral gyrus \\ \hline
\end{tabular}
\caption{Sources in the ``sensory'' time window (50-200 ms) for target stimuli, ranked by normalized dipole moment amplitude (five-direction analysis, active visual oddball task, one participant). Brodmann areas are indicated in parentheses, where applicable.}
\label{S5}
\end{table}

\begin{table}[h!]
\centering
\footnotesize
\begin{tabular}{|c|c|p{0.55\linewidth}|}
\hline
\textbf{Source No.} & \textbf{MNI coordinates} & \textbf{Approx. Anatomical Location} \\ \hline
1 & -40.66, -35.91, -21.44 & Left fusiform gyrus (BA 20) \\ \hline
2 & 9.80, 42.34, 42.55 & Right medial part of the superior frontal gyrus \\ \hline
3 & -57.58, 2.06, 26.47 & Left precentral gyrus (BA 9) \\ \hline
4 & -29.77, -22.30, 8.83 & Left insula \\ \hline
5 & -62.22, -22.50, -21.26 & Left inferior temporal gyrus (BA 20) \\ \hline
6 & -46.65, 23.00, 32.20 & Left middle frontal gyrus \\ \hline
7 & 11.48, 29.98, 48.84 & Right medial superior frontal gyrus \\ \hline
8 & -27.21, 10.96, 54.79 & Left middle frontal gyrus \\ \hline
9 & -35.72, 31.23, 43.89 & Left middle frontal gyrus \\ \hline
10 & -21.03, 2.45, 50.57 & Left middle frontal gyrus \\ \hline
11 & -44.77, 33.17, 23.47 & Left inferior frontal gyrus, triangular part \\ \hline
12 & -37.68, 11.40, 26.22 & Left inferior frontal gyrus, triangular part \\ \hline
13 & -27.30, -29.60, 46.49 & Left superior parietal lobule \\ \hline
14 & -7.76, -56.00, 48.16 & Left precuneus (BA 7) \\ \hline
\end{tabular}
\caption{Sources in the ``cognitive'' time window (350-600 ms) for target stimuli, ranked by normalized dipole moment amplitude (six-direction analysis, active visual oddball task, one participant). Brodmann areas are indicated in parentheses, where applicable.}
\label{C6}
\end{table}
\pagebreak

\begin{table}[H]
\centering
\footnotesize
\begin{tabular}{|c|c|p{0.55\linewidth}|}
\hline
\textbf{Source No.} & \textbf{MNI coordinates} & \textbf{Approx. Anatomical Location} \\ \hline
1 & -40.66, -35.91, -21.44 & Left fusiform gyrus (BA 20) \\ \hline
2 & 9.80, 42.34, 42.55 & Right superior frontal gyrus, medial segment \\ \hline
3 & -52.49, -47.95, 48.57 & Left inferior parietal lobule (BA 40) \\ \hline
4 & 16.97, 23.94, 44.52 & Right superior frontal gyrus \\ \hline
5 & 21.06, 33.77, 52.16 & Right superior frontal gyrus \\ \hline
6 & -52.87, 5.86, 29.87 & Left precentral gyrus \\ \hline
7 & -36.06, -37.07, 57.37 & Left postcentral gyrus \\ \hline
8 & 50.65, -15.60, 54.22 & Right postcentral gyrus (BA 3) \\ \hline
9 & 28.32, 10.98, 55.73 & Right middle frontal gyrus \\ \hline
10 & 40.20, -22.29, -32.73 & Right fusiform gyrus (BA 20) \\ \hline
11 & -57.58, 2.06, 26.47 & Left inferior frontal gyrus (BA 9) \\ \hline
12 & -27.95, 31.25, 32.73 & Left middle frontal gyrus \\ \hline
13 & -29.77, -22.30, 8.83 & Left insula \\ \hline
14 & -62.22, -22.50, -21.26 & Left inferior temporal gyrus (BA 20) \\ \hline
15 & -46.65, 23.00, 32.20 & Left middle frontal gyrus \\ \hline
16 & -38.39, -22.79, 58.11 & Left precentral gyrus (BA 4) \\ \hline
17 & 11.48, 29.98, 48.84 & Right superior medial frontal gyrus \\ \hline
18 & 30.45, 48.67, 26.95 & Right middle frontal gyrus (BA 10) \\ \hline
19 & -39.76, 14.97, 47.67 & Left middle frontal gyrus \\ \hline
20 & -18.67, -8.28, -28.54 & Left parahippocampal gyrus (BA 35) \\ \hline
21 & -44.10, -0.71, -38.36 & Left inferior temporal gyrus \\ \hline
22 & 58.08, -17.48, -34.85 & Right inferior temporal gyrus \\ \hline
23 & -27.21, 10.96, 54.79 & Left middle frontal gyrus \\ \hline
24 & 63.74, -23.36, -7.93 & Right middle temporal gyrus (BA 21) \\ \hline
25 & 56.15, 1.55, -29.62 & Right middle temporal gyrus (BA 21) \\ \hline
26 & -44.29, 39.46, -10.05 & Left inferior frontal gyrus, orbital part (BA 47) \\ \hline
27 & -35.72, 31.23, 43.89 & Left middle frontal gyrus \\ \hline
28 & -21.03, 2.45, 50.57 & Left middle frontal gyrus \\ \hline
29 & 51.51, -67.62, 32.03 & Right angular gyrus (BA 39) \\ \hline
30 & -23.26, 64.98, -4.29 & Left superior frontal gyrus, orbital part (BA 10) \\ \hline
31 & 9.28, -12.52, 65.98 & Right supplementary motor area (BA 10) \\ \hline
32 & 15.36, 52.44, -13.54 & Right superior frontal gyrus, orbital part \\ \hline
33 & -44.77, 33.17, 23.47 & Left inferior frontal gyrus, triangular part \\ \hline
34 & -22.03, -54.91, 48.93 & Left superior parietal lobule (BA 7) \\ \hline
35 & -37.68, 11.40, 26.22 & Left inferior frontal gyrus, triangular part \\ \hline
36 & -27.30, -29.60, 46.49 & Left superior parietal lobule \\ \hline
37 & 59.58, 2.58, -24.20 & Right middle temporal gyrus (BA 21) \\ \hline
38 & -31.13, 52.93, 24.34 & Left middle frontal gyrus \\ \hline
40 & -7.76, -56.00, 48.16 & Left precuneus (BA 7) \\ \hline
41 & 34.22, -15.46, 61.96 & Right precentral gyrus \\ \hline
\end{tabular}
\caption{Sources in the ``cognitive'' time window (350-600 ms) for target stimuli, ranked by normalized dipole moment amplitude (ten-direction analysis, active visual oddball task, one participant). Brodmann areas are indicated in parentheses, where applicable.}
\label{C10}
\end{table}
\pagebreak

\begin{figure}[h!]
\centering
\includegraphics[width=0.8\linewidth]{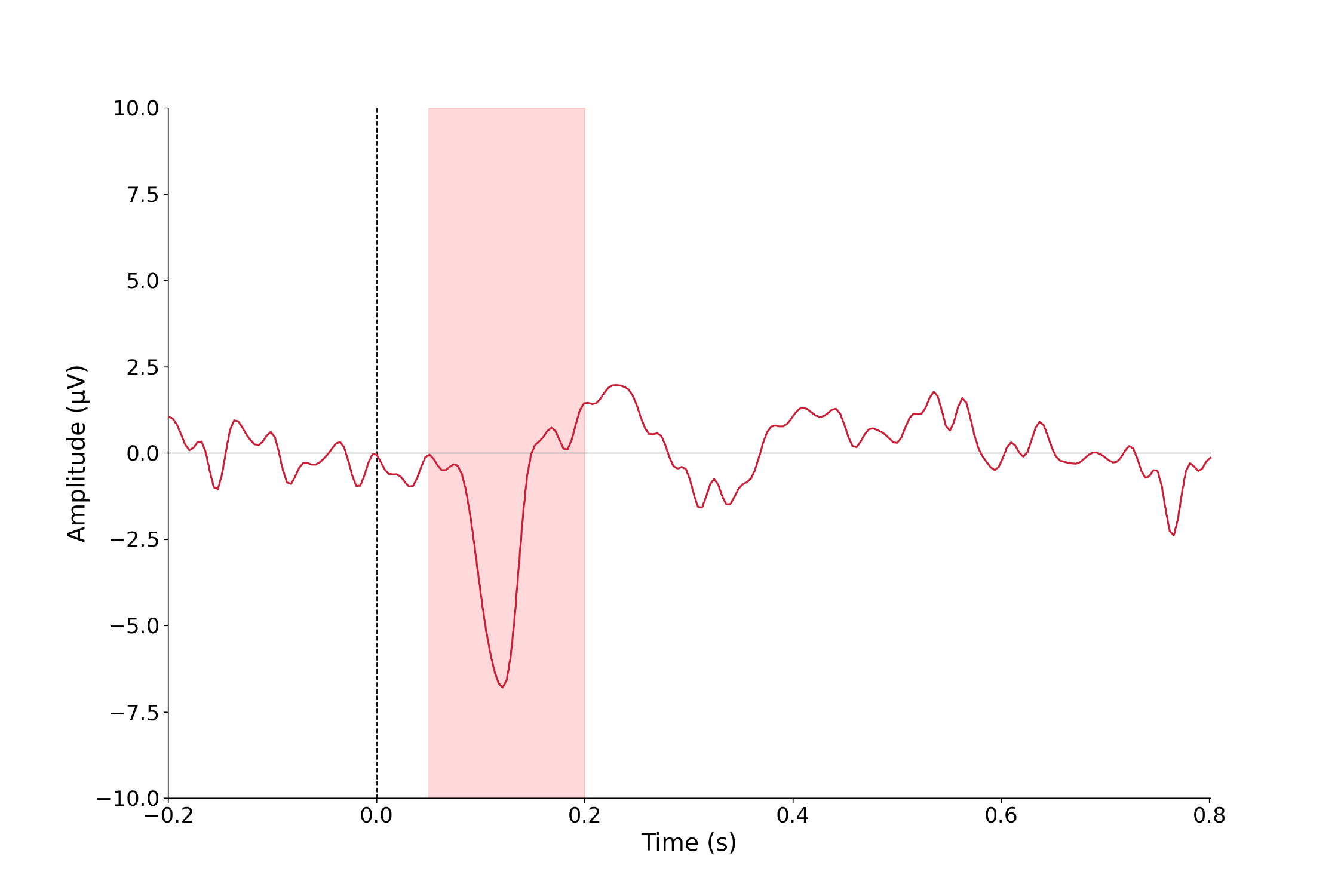}
\caption{The ``sensory'' $[50, 200]$ ms interval at the Oz channel. For visualization purposes, the data were referenced to linked TP9 and TP10 electrodes.} \label{Oz}
\end{figure}

\begin{figure}[h!]
\centering
\includegraphics[width=0.8\linewidth]{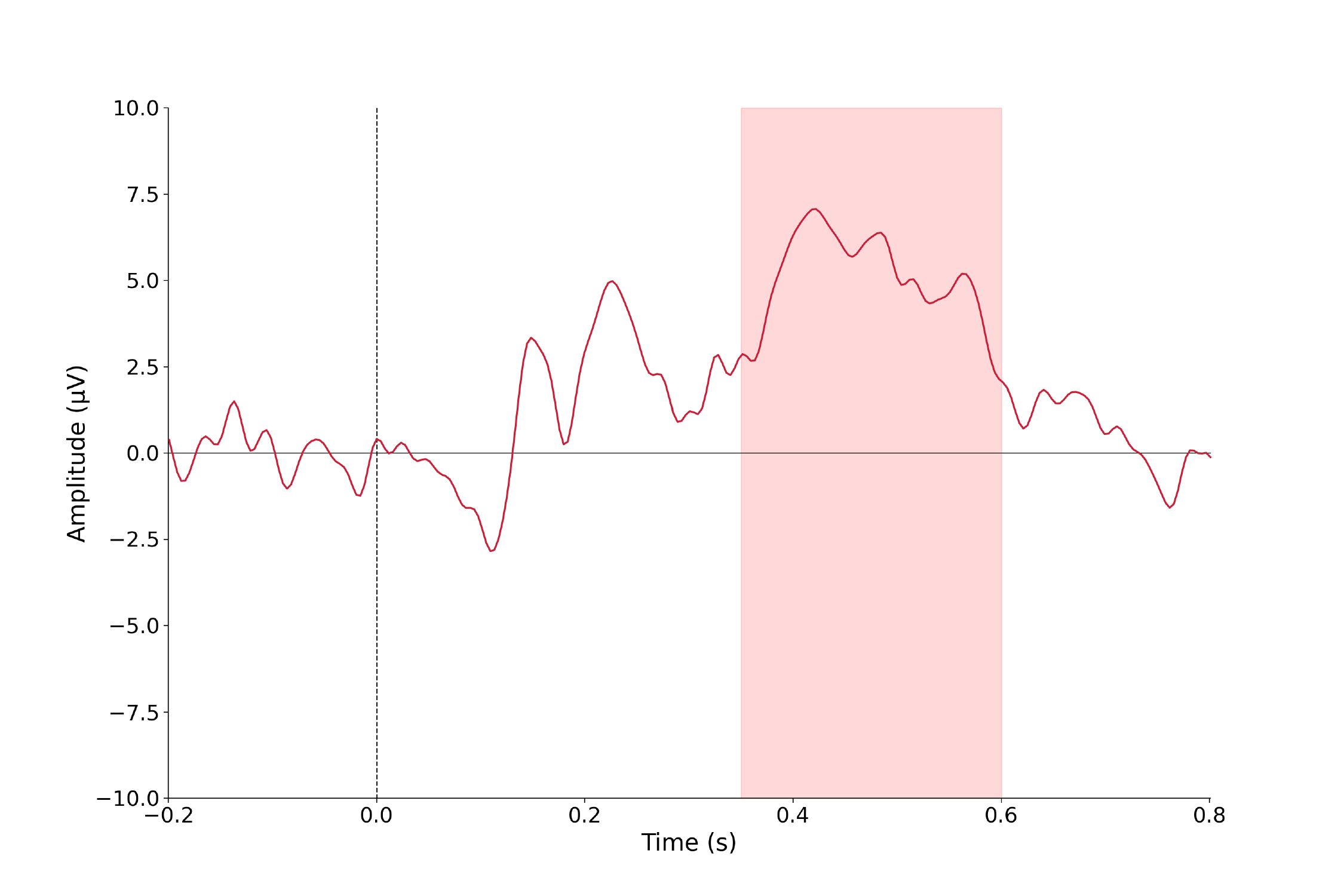}
\caption{The ``cognitive'' $[350, 600]$ ms interval at the Pz channel. For visualization purposes, the data were referenced to linked TP9 and TP10 electrodes.} \label{Pz}
\end{figure}

\section*{Declaration of generative AI and AI-assisted technologies in the manuscript preparation process}
During the preparation of this manuscript, the authors used generative AI to improve the clarity and quality of the English, as none of the co-authors is a native speaker. AI-assisted tools were also used to support the verification of the correctness and efficiency of the code implementing the proposed algorithm. All content generated or refined using these tools was carefully reviewed and edited by the authors, who take full responsibility for the final version of the manuscript and the accompanying code.

\pagebreak
\bibliographystyle{elsarticle-num} 
\bibliography{references}

\end{document}